\documentclass[10pt,a4paper,twocolumn]{article}
\usepackage{setspace}

\usepackage[utf8]{inputenc}
\usepackage[english]{babel}
\usepackage{mathrsfs, amsmath}
\usepackage{amsfonts}
\usepackage{amssymb}
\usepackage{graphicx}
\usepackage[left=2cm,right=2cm,top=2cm,bottom=2cm]{geometry}

\usepackage[toc,page]{appendix}
\usepackage{url}
\usepackage{subcaption}
\usepackage{authblk}
\usepackage{cleveref}
\usepackage[makeroom]{cancel}

\usepackage{bm}
\usepackage{xcolor}

\author[1,2,$\dagger$]{Blas Kolic}
\author[1,2,3,4]{Juan Sabuco}
\author[1,2,5]{J. Doyne Farmer}

\affil[1]{{Institute for New Economic Thinking at the Oxford Martin School, University of Oxford}}
\affil[2]{{Mathematical Institute, University of Oxford}}
\affil[3]{{School of Geography and the Environment, University of Oxford, Oxford, OX1 3QY, UK}}
\affil[4]{{Smith School of Enterprise and the Environment, University of Oxford}}
\affil[5]{{Santa Fe Institute, Santa Fe}}
\affil[$\dagger$]{\textit{blas.kolic@maths.ox.ac.uk}}
 
\title{Estimating initial conditions for dynamical systems with incomplete information}

\newcommand{\xk}[1]{\bm{x}_{#1}}

\newcommand{\model}[1]{\bm{f}(#1)} 
\newcommand{\modelk}[2]{\bm{f}^{ \{m #2\} }(#1)} 

\newcommand{\y}{\bm{y}}
\newcommand{\nse}{\text{NSE}}

\begin{document}
\maketitle

\begin{abstract}
In this paper, we study the problem of inferring the latent initial conditions of a dynamical system under incomplete information, i.e. we assume we observe aggregate statistics of the system rather than its state variables directly. Studying several model systems, we infer the microstates that best reproduce an observed time series when the observations are sparse, noisy and aggregated under a (possibly) nonlinear observation operator. This is done by minimizing the least-squares distance between the observed time series and a model-simulated time series using gradient-based methods. We validate this method for the Lorenz and Mackey-Glass systems by making out-of-sample predictions. Finally, we analyze the predicting power of our method as a function of the number of observations available. We find a critical transition for the Mackey-Glass system, beyond which it can be initialized with arbitrary precision.

\end{abstract}



\section{Introduction}
\label{sec:intro}

We often model empirical processes as complex systems: they are composed of different agents or \textit{microstates}, that interact through simple rules but evolve in nontrivial ways. A lot of research has taken a complex systems approach, where authors have constructed models from the bottom up, but where empirical observations exhibit incomplete information of the system. Applications range from Earth's global climate \cite{Climate1999}, to urban dynamics \cite{ward2016dynamic}, to the human brain \cite{Brain2009}, to financial markets \cite{Markets2012}, to transportation networks \cite{Networks2018}, and marine fisheries \cite{glaser2014complex}. 
However, data are usually available as an aggregate statistic of the microstates because it is difficult to measure microstates directly. Thus, it is an ongoing challenge to develop methods to recover the latent microstates from aggregate and incomplete observations \cite{grazzini2017bayesian}. 

Estimating and forecasting complex systems accurately depends on 1) the inherent complexity of the system, which depends on things like the systems' state space dimension, its Lyapunov exponents, and its attractors,  2) the sparsity and quality of the data, and 3) our ability to model them. In low-dimensional systems with high-quality data, as shown by Packard \textit{et al.} \cite{Farmer1980} and Takens \cite{takens1981detecting}, state-space reconstruction techniques can be applied to partial information to create a representation as a dynamical system. By reconstructing an attractor from data, we can make accurate predictions by choosing its closest points to the current state of the system and extrapolating them \cite{farmer1987predicting}. These techniques work well even without any modelling \cite{sugihara2015free}. 

In high-dimensional systems, it is typically not possible to use time series models.  It is nonetheless often possible to use a theoretical model, if one can only measure the initial conditions that the model requires and compare them to the data.   The task of estimating initial conditions that match observations is known as \textit{initialization}. Similar to the state-space reconstruction techniques, we require to know the evolution function of the underlying dynamics of the system, or at least some approximation of it \cite{farmer1991noisered}. In particular, if the observations available are an aggregate of the dynamical system, then the process is known as \textit{microstate initialization}, or latent state initialization. To do this we need to know how the microstate is aggregated, in addition to having a model of its dynamics.

In the fields of meteorology and numerical weather prediction (NWP), researchers have developed a framework for estimating the latent states of a system \cite{Pires1996, Carrassi2018}, where a large number of observed states are available. This framework is known as data assimilation, and, although it is well justified from the Bayesian perspective, it incorporates several heuristics pertinent to the weather prediction field. For instance, modelers often incorporate Gaussian priors in the cost functionals involved with known statistics about the latent states \cite{navon2009data}. Ideas from data assimilation have already permeated outside the NWP community, such as in urban dynamics \cite{ward2016dynamic}. Typically, data assimilation methods operate in a sequential matter: the microstate at time $t$ gets nudged (or corrected) so it optimally approaches the empirical observation at $t$. Then, the modeler simulates the nudged microstate from time $t$ to time $t + 1$, where the microstate is again nudged.     
Therefore, the resulting sequence of microstates is not a solution of the underlying dynamical system --the microstate was constantly altered throughout its trajectory.
We, in contrast, seek to analyze real trajectories of the latent states, so this sequential approach does not suit our goal.

The initialization process is an optimization problem where we minimize a cost function that depends on some notion of distance between the observed and the model-generated data. Hence, in high-dimensional systems, it is essential to develop efficient algorithms to find initial conditions with high precision. Research has been done around the parameter estimation of stochastic dynamical systems in low-dimensional parameter spaces \cite{grazzini2017bayesian, platt2020comparison}. Among other alternatives, gradient descent has proven to be superior in strongly nonlinear models \cite{evensen1997advanced, Judd2008}, but the main drawback is that it can get stuck in local minima. Other methods, like genetic algorithms, simulated annealing, and other meta-heuristics algorithms \cite{yang2010metaheuristics} usually find the global minimum in low dimensions, but they are likely to diverge in high-dimensional systems. We provide a thorough comparison of the state-of-the-art gradient descent methods \cite{ruder2016overview} and discuss their performance in the context of microstate initialization. 

Here we propose a gradient-based method for initializing chaotic systems where only aggregate observations of short observation windows are available. Using a combination of numerical simulations and analytical arguments, we study the conditions in which certain systems may be initialized with arbitrary precision. We explore the performance of several gradient descent algorithms and the effects of observational noise on the accuracy and convergence of the initialization process. Furthermore, we quantify the accuracy of our method numerically with out-of-sample forecasts. Under this framework, we offer a better understanding of what information the observations provide about a system's underlying dynamics, and, additionally, lay out the connections of our method to those in the data assimilation literature. 

The remainder of this paper is laid out as follows. In section \ref{sec:methods}, we describe the initialization problem under the framework of
dynamical systems and develop our initialization methodology accordingly. In section \ref{sec:results}, we test our method in two systems, namely the Lorenz and Mackey-Glass systems. Finally, in section \ref{sec:conclusion}, we discuss our results and suggest further research directions. 

\section{Methods}
\label{sec:methods}

\subsection{Problem setup}
\label{sec:description} 

A wide class of dynamical systems can be formulated as follows 
\begin{equation}
    \bm{x}(t + \Delta t) = \bm{f}( \bm{x}(t); \bm{\theta}) + \bm{\xi}(t)  \ ,
    \label{eq:ds}
\end{equation}
where $\bm{x}(t) \in 
\mathcal{X} 
\subset \mathbb{R}^{N_x}$ is the ($N_x$-dimensional) state of the system at time $t \in \mathbb{R}$ living in some manifold $\mathcal{X}$, $\bm{\xi}(t) \in \mathbb{R}^{N_x}$ are random variables that model the intrinsic noise of the system and $\bm{\theta} \in \mathbb{R}^d$ is a vector of parameters. 
We will hereafter refer to the state $\bm{x}$ as the \textit{microstate} to emphasize that $\bm{x}$ is not directly observable and is potentially high dimensional.
If $\Delta t > 0$ is finite, the model $\bm{f}: \mathcal{X} \times \mathbb{R}^d \to \mathcal{X}$ is a discrete mapping that takes the microstate from time $t$ to time $t + \Delta t$. If $\Delta t$ is infinitesimal, the system (\ref{eq:ds}) defines a continuous-time dynamical system. 
Here we assume that the dynamical system $\bm{f}$ is deterministic and perfectly specified, so that Eq. (\ref{eq:ds}) depends on neither $\bm{\xi}$ nor $\bm{\theta}$.

We are interested in using information from the past to estimate the microstate $\bm{x}(t_0)$ of a dynamical system $\bm{f}$ at the present time $t_0$.  We assume we know $\bm{f}$, but cannot directly observe its microstate $\bm{x}(t)$. Instead, we are only able to observe a sequence of observations $\y = (y_{-T}, \dots, y_0)$ measured at times $\{t_{-T}, \dots, t_0 \} $, where $y_k \in \mathbb{R}^{N_y}$ and $N_y < N_x$, i.e., we are interested in measurements that lose information of the microstate by reducing its dimension from $N_x$ to $N_y$. Moreover, these observations can be noisy. Thus, we can relate the observations to the dynamical system as
\begin{equation}
    y_{k} =  \mathcal{H} \big( \bm{x}(t_k) \big) + \epsilon_{k} \ .
    \label{eq:observations}
\end{equation}
where $\mathcal{H}: \mathcal{X} \subset \mathbb{R}^{N_x} \to \mathbb{R}^{N_y}$ is a known (possibly) nonlinear \textit{observation operator} and $\epsilon_k$ accounts for \textit{observational noise}. Although it is not  necessary for the following discussions, we will assume for simplicity that $N_y = 1$ and, that the noise has a known variance $\sigma_y$. Our choice of indices from $-T$ to $0$ for $y_k$ reflects our interest in the \textit{present time} $t_0$.

We assume that the observations are sampled at a uniform rate as follows
\begin{equation*}
    \Delta t_k = t_{k+1} - t_k = m \Delta t \ ,    
\end{equation*}
where the non-zero integer $m$ is the \textit{sampling interval}\footnote{The initialization procedure we introduce in the following section works with irregularly sampled observations as well.}. We treat $m$ as a parameter that controls how often we sample observations from the system.
Thus, we can recast Eq. (\ref{eq:ds}) as a discrete-time mapping
\begin{equation}
    \bm{x}_{k+1} = \bm{f}^{ \{m \} }( \bm{x}_k ) = \bm{f}^{ \{mk \} }( \bm{x}_{-T}) \ , 
    \label{eq:dynamical_map}
\end{equation}
where $\bm{f}^{ \{ i \} }(\cdot)$ is the  composition of $\bm{f}$ with itself $i$-times, $\bm{x}_{-T} = \bm{x}(t_{-T})$ is the latent microstate at the time of the first observation, and $\bm{f}^{ \{k \} }$ is the self-iteration of $f$ 
by $k$ times. Note that the whole evolution of the microstates is determined by $\bm{x}_{-T}$, so the \textit{microstate initialization problem} is that of obtaining the best approximation to $\bm{x}_0$ given the observations from the \textit{assimilation time} $t_{-T}$ up to present time $t_0$. These observations, alongside with the model (\ref{eq:dynamical_map}) and the observation operator (\ref{eq:observations}), form a nonlinear system of equations with $N_x$ variables and $T$ equations corresponding to the dimension of the microstate space and the number of observations.

\subsection{Initialization procedure}
\label{sec:minimization_procedure}

As we stated previously, the microstate initialization problem is that of obtaining the best representation of the present-time microstate $\bm{x}_0$ given the history of observations $\y$ under the dynamical system $\bm{f}$. We define what \textit{best} means in what follows. We make the code of the initialization procedure freely available at \url{https://github.com/blas-ko/LatentStateInitialization}.

\subsubsection{Cost function}\label{sec:cost_func}

Given an estimate $\bm{x}$ of the ground-truth microstate $\bm{x}_{-T}$ at the assimilation time $t_{-T}$, a natural way to quantify the goodness of $\bm{x}$ is by measuring the point-to-point discrepancy between the observed and estimated data. We do so with the following mean-squared \textit{cost function}.\footnote{We may generalize this cost function for non-scalar observations as $$\mathcal{J}(\bm{x}) = \frac{1}{T} \sum_k (\bm{y}_k - \hat{\bm{y}}_k)^\dagger\bm{\Sigma}_y^{-1}(\bm{y}_k - \hat{\bm{y}}_k) \ ,$$ where $\bm{\Sigma}_y$ is the covariance matrix of the observations.}
\begin{equation}
     \mathcal{J}( \bm{x} ) = \frac{1}{T \sigma_y^2} \sum\limits_{k = -T}^{0} \left( y_k - \hat{y}_k  \right)^2 ,
    \label{eq:cost_function}
\end{equation}
where
\begin{equation}
     \hat{y}_k(\bm{x}) = \mathcal{H} \Big( \bm{f}^{ \{ mk \} }( \bm{x} ) \Big)
     \label{eq:???}
\end{equation}
is our estimate of observation $y_k$ given $\bm{x}$, and $\sigma_y$, the variance of the observed data, is a normalization constant that is arbitrary but will be of practical use later. We call the state $\hat{\bm{x}}_{-T}$ that minimizes $\mathcal{J}$ the \textit{assimilated microstate} and the present-time state $\hat{\bm{x}}_0 = \modelk{\hat{\bm{x}}_{-T}}{T}$ the \textit{initialized microstate}.  Our goal is to find an initialized microstate $\hat{\bm{x}}_0$ that is a good representation of the ground-truth microstate $\bm{x}_0$.

The cost function $\mathcal{J}$ is known as a \textit{filter} in the interpolation assimilation community, a \textit{least-squares optimizer} in the optimization and machine learning communities, and \textit{4D-Var} with infinite state uncertainty in the variational assimilation community. From a Bayesian perspective, the cost function $\mathcal{J}$ emerges naturally when we assume the prior $p(\bm{x})$ is uniform and observational noise is Gaussian: the posterior ${p(\bm{x} | \y )}$ is maximal when $\mathcal{J}$ is minimal \cite{Carrassi2018}.

If the observations in the time series $\y$ are noiseless, 
and given that we assume that the model $\bm{f}$ perfectly describes the system,
then $\mathcal{J}$ has a global minimum $\hat{\bm{x}}_{-T}$ for which $\mathcal{J}(\hat{\bm{x}}_{-T}) = 0$. However, even for the noiseless scenario, $\hat{\bm{x}}_{-T}$ is not necessarily unique. Take, for instance, the Lorenz system \cite{lorenz1963deterministic}, which is symmetric around its $x$-axis, and take an observation operator of the form $\mathcal{H}(\bm{x}) = \bm{x}^\dagger \bm{x}$ with $\dagger$ denoting matrix transposition. Note that for any microstate $\bm{x}$ in the axis of symmetry, the condition $-\bm{x}$ will produce the same sequence of noiseless observations, so $\bm{x}$ and $-\bm{x}$ are indistinguishable in terms of $\mathcal{J}$. We will refer to the set of indistinguishable microstates as the \textit{feasible set} of solutions, and we will denote it as $\bm{\Omega}$. 

An additional problem arises when the observations contain noise. In any finite time series, there might exist microstates with a lower value of the cost function than the ground-truth microstate --i.e., where $\mathcal{J}( \hat{\bm{x}}_{-T} ) < \mathcal{J}( \bm{x}_{-T} ) $ for $\hat{\bm{x}}_{-T} \not\in \bm{\Omega}$-- which we will call \textit{dominating microstates}, following \cite{Judd2007}. Judd \textit{et al.} \cite{Judd2007} suggest that using cost functions that minimize both the variance (as in Eq. (\ref{eq:cost_function})) and the kurtosis will do a better job of identifying the true microstate when the observational noise is Gaussian. While this is a good idea when there are many observations, we consider short time series in this work, which calls for other alternatives. Instead of modifying the cost function $\mathcal{J}$, we  pre-process the data to reduce the probability of finding dominating trajectories outside of the feasible set $\bm{\Omega}$. 

Following \cite{Pires1996}, we distinguish between the \textit{error-free} and the \textit{noise} contributions to the cost function $\mathcal{J}$. Isolating the contributions from the discrepancy between ground-truth and assimilated microstate and the observational noise will let us design a well-suited methodology to deal with noise and dominating trajectories in a separate manner. Recall that, given $\bm{x}_{-T}$, $\mathcal{H}(\bm{x}_k)$ is the error-free observation at time $t_k$, $\epsilon_k$ its associated noise and $\hat{y}_k$ our estimation of $y_k$. Thus, we can decompose Eq. (\ref{eq:cost_function}) into three terms of the form
\begin{align*}
    \mathcal{J}(\bm{x}) = \frac{1}{\sigma_y^2} \Big[& \underbrace{\frac{1}{T}\sum_k (\mathcal{H}(\bm{x}_k) - \hat{y}_k)^2}_{  \text{noise-free cost}} \\
    +& \underbrace{\frac{1}{T}\sum_k \epsilon_k^2}_{\text{obseravation noise}} \\
    -& \underbrace{\frac{2}{T}\sum_k (\mathcal{H}(\bm{x}_k) - \hat{y}_k)\epsilon_k}_{\text{error-noise covariates}} \Big]\ 
\end{align*}
The first term on the RHS is the noise-free cost, $\mathcal{J}_{\text{free}}(\bm{x})$, that we would obtain in the absence of observational error. The second term is the average contribution of the square of the noise, which converges to the noise variance, $\sigma_n^2$, for large $T$. Finally, the third term captures how the noise and the noise-free cost vary together, which goes to $0$ for large $T$ because we assume the observational noise and the system dynamics are uncorrelated. Considering a large number of observations, we can thus approximate $\mathcal{J}$ as 
\begin{equation}
    \mathcal{J}(\bm{x}) \approx \mathcal{J}_{\text{free}}(\bm{x}) + \frac{\sigma_n^2}{\sigma_y^2} \ ,
    \label{eq:expected_cost}
\end{equation}
which shows the expected behavior of Eq. (\ref{eq:cost_function}) in the presence of noise.

Now, note that if we take $\hat{y}_k = \mathbb{E}[\bm{y}]$ for all $k$ \cite{farmer1987predicting}, which is the best constant predictor for the observed time series, then $\mathcal{J}_{\text{free}} = 1$ and, therefore 
\begin{equation*}
    \mathcal{J}(\bm{x} : \hat{y}_k = \mathbb{E}[\bm{y}] \ \forall k) = \mathcal{J}_{\text{const}} := 1 + \sigma_n^2 / \sigma_y^2.
\end{equation*}
Naturally, we want to find an assimilated microstate $\hat{\bm{x}}_{-T}$ that performs better than a constant predictor, so that $\mathcal{J}(\hat{\bm{x}}_{-T}) \leq \mathcal{J}_{\text{const}}$. This motivates us to consider microstates $\bm{x}$ such as
\begin{equation}
    \lbrace \bm{x}: \mathcal{J}(\bm{x}) \leq \alpha + \frac{\sigma_n^2}{\sigma_y^2}\beta \rbrace \ ,
\end{equation}
where the parameter $\alpha \in (0,1]$ accounts for the error-free tolerance about the global minimum while $\beta \in (0,1]$ accounts for the noise tolerance. We will explore these parameters in the following sections.

The above discussion suggests that we should pay especial attention on handling dominating trajectories, which might be present by either the presence of observational noise or by microstates that live far away from the ground-truth but that have low cost function values. Thus, we propose the following three stages to initialize the microstate from aggregate observations:
1) a \textit{preprocessing} stage in which we reduce the noise of the observations, 2) a \textit{bounding} stage in which we limit the region of the microstate space in which we search for an optimal solution, and 3) a \textit{refinement} stage in which we minimize the cost function (\ref{eq:cost_function}) in a small search space and estimate the optimal microstate given the observations. 

\subsubsection{Preprocess: noise reduction}\label{sec:preprocess}

First, we preprocess the observed time series to reduce the observational noise and thus lower the probability of obtaining dominating microstates. Casdagli \textit{et al.} \cite{Casdagli1991} showed that in the presence of observational noise, the distribution of local minima gets increasingly complex with increasing levels of noise, especially when the dynamics is chaotic. We handle time series with only a handful of data points (in the order of $100$ data points or less), so reducing noise by orbit shadowing \cite{farmer1991noisered} or large window impulse response filters \cite{guinon2007moving} are not suitable options. Instead, we find that the best way to reduce the variance of the noise is using a low-pass moving average (LPMA) filter,
\begin{equation}
    z_k := \begin{cases} 
      \frac{1}{2} y_{k} + \frac{1}{2} y_{k+1} & k = -T \\
      \frac{1}{2} y_{k} + \frac{1}{2} y_{k-1} & k = 0 \\
      \frac{1}{2} y_k + \frac{1}{4} \left( y_{k-1} + y_{k+1} \right) & \text{otherwise} \\
   \end{cases} \ ,
    \label{eq:noise_filter} 
\end{equation}
where $z_k$ is the filtered data point at time $t_k$. 
We control the amount of noise reduction by repeatedly feeding the signal back into the LPMA filter of Eq. (\ref{eq:noise_filter}). Feeding the signal back into the filter $q$ times, hereafter denoted as $z_k^q$, is equivalent to increasing the filtering window from three to $2q+1$ points, making the filtered signal smoother. 

We expect for the resulting variance, $(\sigma_n^q)^2$, to be lower than the original noise variance $\sigma_n^2$. Thus, following \cite{grassberger1993noise} and assuming we can rewrite the filtered signal in terms of the microstates as
$z_k^q = \mathcal{H}(\bm{x}_k) + \epsilon_k^q$, we can measure the performance of the LPMA filter with the increase of the signal-to-noise ratio
\begin{equation}
    r_0 = \sqrt{  \frac{\sum_k (\epsilon_k)^2 }{ \sum_k (\epsilon_k^q)^2 } } \approx \frac{\sigma_n}{\sigma_n^q} ,
    \label{eq:snr_increase}
\end{equation}
where $r_0 > 1$ whenever $\sigma_n > \sigma_n^q$.

The resulting noise distribution of the filtered signal converges to a zero-mean Gaussian distribution if we have either many data samples or set $q$ big enough. However, if $q$ is too big, we may filter parts of the dynamics and mix them into noise, resulting in exotic noise distributions (see Fig. \ref{fig:error_distributions} in Appendix \ref{app:interesting_figs} for examples). What \textit{big enough}, \textit{many samples} and \textit{too big} mean depend heavily on the dynamical system and the noise distribution, although we stress that the LPMA filter works optimally when the noise distribution has higher frequency spectrum on average than that of the dynamics of the system (see \cite{natrella2010nist}, Chapter 6.4.3.1.).

\subsubsection{Bound: exploring the attractor}\label{sec:bound}

After the preprocessing stage, the next step is to bound the search space. Under no constraints in the cost function, the initialization procedure consists on searching through the whole microstate space $\mathcal{X}$ for a set of microstates that minimize Eq. (\ref{eq:cost_function}), with no prior preference on where to start the search from. However, many real world systems are dissipative, meaning that their dynamics relax into an attractor, i.e. a subset manifold $\mathcal{M} \subset \mathcal{X}$ of the microstate space. Thus, it is safe to assume that the observations $\y$ derive from a sequence of microstates that live in or near $\mathcal{M}$, and, by the properties of dissipative systems, any microstate $\bm{x}$ in the basin of attraction will eventually visit every point in $\mathcal{M}$ \cite{anishchenko2002peculiarities}. This means that, if we wait for long enough, then any point in the basin of attraction will get arbitrarily close to the ground truth microstate.

The bounding stage consists of exploiting the dissipative nature of real world systems and letting any arbitrary estimate of the microstate explore the basin of attraction until it \textit{roughly} approaches the ground truth microstate. To be more precise, we say that the microstate $\bm{x}_{-T}^R \in \mathcal{X}$ \textit{roughly approaches} $\bm{x}_{-T}$ if    
\begin{equation}
    \mathcal{J}(\bm{x}_{-T}^R) \leq \delta_R := \alpha_R + \frac{\sigma_n^2}{\sigma_y^2}\beta_R, 
    \label{eq:deltaR}
\end{equation}
for some \textit{rough threshold} $0 \ll \delta_R < 1$. Thus, we let an arbitrary microstate evolve according to the model $\bm{f}$ until Eq. (\ref{eq:deltaR}) is satisfied. 

At this stage, we want to obtain solutions with a cost value of the order of the unfiltered noise level $\sigma_n^2/\sigma_y^2$, so that $\bm{x}_{-T}^R$ is either near to the feasible set $\bm{\Omega}$ or to any of the dominating microstates driven by the noise. To achieve this, it suffices to set $\beta_R \sim \mathcal{O}(1)$ and $\alpha_R < \sigma_n^2/\sigma_y^2 \leq 1$. 

We note that whenever the time series $\bm{y}$ is noiseless, the bounding stage is only driven by $\alpha_R$, the error-free tolerance of the points situated at the global minima of the cost function, which are exactly those in the feasible set $\bm{\Omega}$. Thus, our choice of $\alpha_R$ leverages how closely we approach to $\bm{\Omega}$. If we set $\alpha_R$ too close to $0$, we would impose for $\hat{\bm{x}}_{-T}$ to lay near $\bm{\Omega}$; however, it would take too long simulation times to satisfy Eq. (\ref{eq:deltaR}) for this approach to be practical. The idea, ultimately, is to set the lowest $\delta_R$ possible such that the time to satisfy Eq. (\ref{eq:deltaR}) is \textit{short}.

\subsubsection{Refine: cost minimization}\label{sec:refine}

The final step of the initialization procedure is refining $\bm{x}_{-T}^R$. By this point, we expect that $\bm{x}_{-T}^R$, our estimate of $\bm{x}_{-T}$, has bypassed most of the high-valued local minima of the cost function landscape. Additionally, we preprocessed the observations $\bm{y}$ to reduce their observational noise, but we have not fully exploited such preprocessing yet. By reducing the variance of the observational noise, we lower the number of dominating trajectories of the cost function --i.e., trajectories for $\bm{x} \not\in \bm{\Omega}$ such that $\mathcal{J}(\bm{x}) < \mathcal{J}(\bm{x}_{-T})$. Thus, starting from $\bm{x}_{-T}^R$, we can minimize $\mathcal{J}$ using any optimization scheme until
\begin{equation}
    \mathcal{J}(\hat{\bm{x}}_{-T}) \leq \delta_r := \alpha_r + \frac{\sigma_n^2}{\sigma_y^2}\beta_r
    \label{eq:deltar}
\end{equation}
for some \textit{refinement threshold} $0 < \delta_r \ll 1$, with $\hat{\bm{x}}_{-T}$ the assimilated microstate of the system. We then define $\hat{\bm{x}}_0 = \modelk{ \hat{\bm{x}}_{-T} }{T}$ to be the initialized microstate, hoping that $\hat{\bm{x}}_0$ is a good representation of $\bm{x}_0$.

We want for $\hat{\bm{x}}_{-T}$ to have \textit{the lowest cost possible} at this stage. In terms of the error-free tolerance, we look for $\alpha_r \ll \alpha_R < \sigma_n^2/\sigma_y^2 \leq 1$, but the actual magnitude of $\alpha_r$ is left to the modeller to choose. Regarding the contribution of the noise, recall from Eq. (\ref{eq:snr_increase}) that $r_0 \approx \sigma_n / \sigma_n^q > 1$, so the lowest expected cost we can get is $\sigma_n^2/\sigma_y^2 r_0^{-2}$ (see Eq. (\ref{eq:expected_cost})). Thus, we define the \textit{refinement bound} of our initialization procedure as that of setting $\alpha_r \ll \alpha_R$ and $\beta_r \sim \mathcal{O}(r_0^{-2})$.

For our optimization scheme we explore a plethora of the most successful \textit{gradient-based algorithms} in the literature \cite{ruder2016overview}. These algorithms include Stochastic Gradient Descent \cite{robbins1951stochastic}, Momentum Descent \cite{polyak1964some}, Nesterov \cite{nesterov1983method}, Adagrad \cite{duchi2011adaptive}, Adadelta \cite{zeiler2012adadelta}, Rmsprop \cite{tieleman2012lecture}, Adam \cite{kingma2015adam}, and AdamX \cite{reddi2019convergence} and YamAdam \cite{yamada2018yamadam}. 

In all cases we set the hyper-parameters to be those given in the literature. We then compute the gradient of $\mathcal{J}$ using centered finite differences of step size $\sqrt{\varepsilon_M} \approx 1.5 \times 10^{-8}$, where $\varepsilon_M$ corresponded to double precision arithmetic on our machine. In the absence of observational noise, 
provided that the dynamics is not degenerate and that 
we have sufficient
observations, we expect that
the feasible set $\bm{\Omega}$ collapses to 
the ground-truth microstate only, so at this stage we expect to infer it with high precision from the data.


\subsection{Validation}
\label{sec:validation}

We validate the initialized microstate $\hat{\bm{x}}_0$ by comparing them with the present-time microstate $\bm{x}_0$ and making out-of-sample predictions of the observed time series. Recall that we take the convention in which the observations, $\y = (y_{-T}, \dots, y_{0})$, have non-positive time indexes, so we refer to the observation times $t_k$ for $k < 0$ as \textit{assimilative} while we refer to times for $k \geq 0$ as \textit{predictive}. We measure the \textit{discrepancy} between the real and the simulated observations using both the normalized squared error in the observation space and the normalized error in the model space, i.e.
\begin{align}
    \nse_k^{obs} &= \frac{  ( y_k - \hat{y}_k )^2  }{\sigma_y^2},
    \label{eq:nse_obs} \\
    \nse_k^{mod} &= \frac{1}{N_x} (\bm{x}_k - \hat{\bm{x}}_k)^\dagger\Sigma_{\bm{x}}^{-1}(\bm{x}_k - \hat{\bm{x}}_k) ,
    \label{eq:nse_mod}
\end{align}
where $\sigma_y^2$ is the variance of the data, $\Sigma_{\bm{x}}$ the covariance matrix of the microstates and $(\cdot)^\dagger$ denotes matrix transposition. In general, $\bm{x}_k$ and $\Sigma_{\bm{x}}$ are unknown to the modeller, but we use them to measure the performance of our initializations in the latent space of the microstates.

Note that if we let $T \to \infty$, then $\mathcal{J}(\bm{x})$ converges to $\mathbb{E} [ NSE^{obs}_0 ]$. When $k$ grows, the trajectories $y_k$ and $\hat{y}_k$ diverge exponentially until they lose all memory about their initial conditions ($\bm{x}_0$ and $\hat{\bm{x}_0}$ respectively). Given that such divergence is exponential, we therefore take the \textit{median} of the $NSE$ when comparing the performance over an ensemble of experiments as a better alternative to the mean.

On the same note, another way to validate the inferred microstate is by looking for how long our predictions accurately describe the system. Chaotic systems are by definition sensitive to initial conditions, meaning that microstates that are close to each other diverge exponentially over time.  If these microstates live in a chaotic attractor, the distance between them is bounded by the size of the attractor \cite{Pires1996}. Thus, we can assess the quality of our predictions by measuring how long the real and simulated observations retain memory about each other \cite{farmer1987predicting}. We refer to this limiting time as the \textit{predictability horizon} of the system $k_{max}$, and we define it as the average number of steps before the separation between $y_k$ and $\hat{y}_k$ is greater than the distance between two random points in the attractor of the system. Mathematically,
\begin{equation*}
    k_{max} := \mathbb{E} \Big[ \arg\min\limits_{ k \geq 0} \big\lbrace ( y_k - \hat{y}_k )^2 \geq \mathcal{D}_\mathcal{M} \big\rbrace \Big],
\end{equation*}
where $\mathcal{D}_\mathcal{M}$ is average squared distance between two random points in the attractor $\mathcal{M}$ normalized by the variance of the attractor. More specifically, if $X,Y$ are two i.i.d. random variables such that $X \sim \mu(\mathcal{M})$ with $\mu(\cdot)$ denoting the natural measure, then 
\begin{equation*}
    \mathcal{D}_\mathcal{M} := \frac{\mathbb{E}[ \lVert X - Y \rVert^2 ] }{ Var(X) } = 2.  
\end{equation*}

Thus, given that $\mathbb{E}[ \sigma_y ] = Var(X)$, we can simplify $k_{max}$ into an expression that only depends on Eq. (\ref{eq:nse_obs}) so that
\begin{equation}
    k_{max} = \mathbb{E} \Big[ \arg\min\limits_{ k \geq 0} \big\lbrace \nse_{k} \geq 2 \big\rbrace \Big] \ ,
    \label{eq:k_max}
\end{equation}
which is the formula we use to compute $k_{max}$.

Finally, we benchmark the inferred microstate by comparing $k_{max}$ with a measure of the natural rate of divergence of the dynamics.  As a measure of this, we use the \textit{Lyapunov 10-fold time} $t_\lambda$, which indicates the average time for two neighboring microstates to diverge from each other by one order of magnitude \cite{aurell1996growth}.  This is defined in terms of the inverse of the maximum Lyapunov exponent $\lambda$ as
\begin{equation}
    t_\lambda  = \frac{ \ln{10} }{m \Delta t} \lambda^{-1},
    \label{eq:lyap_time}
\end{equation}
where the factor $\ln{10} /(m \Delta t)$ lets us interpret $t_\lambda$ in the units of the number of observations after which on average the dynamics causes the loss of an order of magnitude of precision. We obtain $\lambda$ numerically using the two-particle method from Benettin \textit{et al.} \cite{benettin1976kolmogorov}.

\subsection{Initialization procedure summary}\label{sec:summary}

We summarize our \textit{microstate initialization procedure} in the following steps: 
\begin{enumerate}
\item \textbf{Preprocess:} Smooth the observed time series using the LMPA filter (see Eq. (\ref{eq:noise_filter})) or any other suitable noise reduction technique. The smoothed signal will have fewer dominating trajectories \cite{Judd2007} and a simpler distribution of local minima than the full noisy signal \cite{Casdagli1991}. 

\item \textbf{Bound:}  Make an arbitrary guess $\bm{x} \in \mathcal{X}$ of the microstate, and let it evolve under the model $\bm{f}$  until the microstate \textit{roughly approaches} the smoothed observations, i.e., until $\mathcal{J}(\bm{x}_{-T}^R) \leq \delta_R$ for $\bm{x}_{-T}^R = \bm{f}^{ \{mR\} }(\bm{x})$ for some $R \geq 0$. (see Eq. (\ref{eq:deltaR})). If several attractors exist, make one arbitrary guess for each of the different basins of attraction in the system. 

\item \textbf{Refine:} Minimize the cost function $\mathcal{J}$ starting from $\mathbf{x}_{-T}^R$ using Adam gradient descent or any other suitable optimization scheme until $\mathcal{J}(\hat{\bm{x}}_{-T}) \leq \delta_r$ (see Eq. (\ref{eq:deltar})) and call $\hat{\bm{x}}_{-T}$ the \textit{assimilated microstate} and $\hat{\bm{x}}_0 = \modelk{\hat{\bm{x}}_{-T}}{T}$ the \textit{initialized microstate}.

\item \textbf{Validate:} Compute the discrepancy between the real and simulated observations (see Eq. (\ref{eq:nse_obs})) and the predictability horizon $k_{max}$ (see Eq. (\ref{eq:k_max})) on out-of-sample predictions of the system to evaluate the quality of the initialized microstate $\hat{\bm{x}}_0$. If possible, benchmark the predictability horizon with the Lyapunov 10-fold time (see Eq. (\ref{eq:lyap_time})) of the system considered.
\end{enumerate}


\section{Results}
\label{sec:results}

In this section, we test the microstate initialization procedure on two paradigmatic chaotic systems: the well-known Lorenz system \cite{lorenz1963deterministic} and the high-dimensional Mackey-Glass system \cite{mackey1977oscillation}. We approximate both systems using numerical integrators (described in each section that follows), and we take the approximated system as the real dynamical system.

In all cases, we sample the ground-truth microstate $\bm{x}_{-T}$ from the attractor of the system considered. We generate observations using the following nonlinear observation operator
\begin{equation}
    \mathcal{H}(\bm{x}) = \sqrt[3]{\sum_{i=1}^{N_x} (\bm{x})_i^3 },
    \label{eq:observation_operator}
\end{equation}
where $(\bm{x})_i$ is the $i$-th component of $\bm{x}$. Note that, while $\mathcal{H}$ is nonlinear, the mappings $x \to x^3$ and $x \to \sqrt[3]{x}$ are bijective, so there exists a diffeomorphism between $\mathcal{H}$ and any non-degenerate linear operator $\bm{H}: \mathcal{X} \subset \mathbb{R}^{N_x} \to \mathbb{R}$. Additionally to this operator, in Appendix \ref{app:observation_operators} we study the quality of our initialization procedure with several observations operators each with different levels of coupling between the microstate components. The predictions in the observation space are almost identical regardless of what observation operator we use. However, given the symmetry about the $x$-axis of the Lorenz system, the operator with the maximum coupling recovers the right $x$ component, but a reflection of the ground truth of the $y$ and $z$ components.

We test our initialization procedure on both noiseless time series and noisy time series. For the noisy series, we take a zero-mean Gaussian noise distribution $\epsilon_k \sim \mathcal{N}(0, \sigma_n^2)$ for all $k$. Here, $\sigma_n$ represents the noise level of the observations, which we take to be $30 \%$ of the standard deviation of the observed data; i.e., $\sigma_n = 0.3\sigma_y$. We always construct the noiseless and noisy time series from the same ground-truth microstate so that all our results are comparable. 

Although our choice for the initial guess is arbitrary, we initialize our method with a microstate that matches the first observation of $\bm{y}$ exactly. This is, we take our initial guess at random from the set $\{ \bm{x} \in \mathcal{X} : \mathcal{H}(\bm{x}) = y_{-T} \}$. This is straightforward to do for any homogeneous function, such as the one of Eq. (\ref{eq:observation_operator}). 

Throughout the results, we use the parameters and system features described in Table \ref{tab:params} unless otherwise stated. However, we evaluate the performance of our method for various choices of the rough parameter $\delta_R$ (see Fig. \ref{fig:effects_of_bounding} in Appendix \ref{app:interesting_figs}) and the optimizers described in section \ref{sec:refine} (See Figs. \ref{fig:lorenz_optimizers_comparison} and \ref{fig:mackey_optimizers_comparison} in Appendix \ref{app:interesting_figs}). Varying $\delta_R$ determines the effect of bounding the search space into finer-grained regions of the attractor, and we find that when observations are noiseless, the lower $\delta_R$ the better the initialization but the longer it takes to meet condition (\ref{eq:deltaR}). For noisy observations, we find that if we set $\delta_R$ very small, the refinement stage yields no improvement over the initialized microstates obtained. In terms of the optimization schemes of the refinement stage, we find that Adadelta \cite{zeiler2012adadelta} and the various flavors of Adam \cite{kingma2015adam, reddi2019convergence, yamada2018yamadam} outperform all other alternatives, with significantly better results than vanilla Gradient Descent.

\begin{table*}[h!]
\centering
\begin{tabular}{l|llll|lllllll}
                      & $\bm{\alpha_R}$ & $\bm{\alpha_r}$ & $\bm{\beta_R}$ & $\bm{\beta_r}$ & $\bm{T}$ & $\bm{N_x}$ & $\bm{m}$  & $\bm{\sigma_n/\sigma_y}$ & $\bm{q}$ & $\bm{r_0}$  & $\bm{t_\lambda}$ 
\\ \hline
\textbf{Lorenz}       & $0.05$          & $10^{-4}$       & $0.5$          & $0.8 r_0^{-2}$ & $50$     & $3$        & $2$       & $0.3$                    & $4$      & $2.02$      &  $127$    \\ 
\textbf{Mackey-Glass} & $0.05$          & $10^{-5}$       & $0.5$          & $0.2 r_0^{-2}$ & $25$     & $50$       & $2$       & $0.3$                    & $5$      & $2.41$      &  $230$    \\ 
\end{tabular}
\caption{\textbf{Parameters and other quantities}. $\alpha_R$, $\alpha_r$, $\beta_R$, and $\beta_r$ are the parameters of our procedure. $T$ is the number of data points in the time series $\bm{y}$, $N_x$ is the dimension of the microstate space, $m$ is the sampling interval between observations, $q$ is the number of times we feed the signal back into the LPMA filter, $r_0$ is the increase in the signal-to-noise ratio, $\sigma_n/\sigma_y$ is the noise level, and $t_\lambda$ is the Lyapunov 10-fold time of the system,.}
\label{tab:params}
\end{table*}

\subsection{Low-dimensional Example: Lorenz System}
\label{sec:lorenz}

We first test our microstate initialization procedure on the Lorenz system \cite{lorenz1963deterministic}, described by the following differential equations
\begin{align}
\begin{aligned}
	 \dot{x} &= \sigma ( y -x ) \ ,  \\
	 \dot{y} &= x(\rho - z) - y \ , \\	
	 \dot{z} &= xy - \beta z \ ,	 
\end{aligned}
\label{eq:lorenz}
\end{align}
where $\bm{x}_k = \big( x(t_k), y(t_k), z(t_k) \big)$ is the microstate of the system at time $t_k$. The dynamics exhibit chaotic behavior for the parameters $\sigma = 28$, $\rho = 10$, and $\beta = 8/3$. We solve the system using a $4$-th order Runge-Kutta integrator with a fixed step size of $0.01$ units. Under these settings, the system's Lyapunov 10-fold time is $ t_\lambda = 127$ samples while its attractor has a Kaplan-Yorke dimension of $N_{\mathcal{M}} = 2.06$.

\begin{figure}[h!]
\centering
	\includegraphics[width=0.49\textwidth]{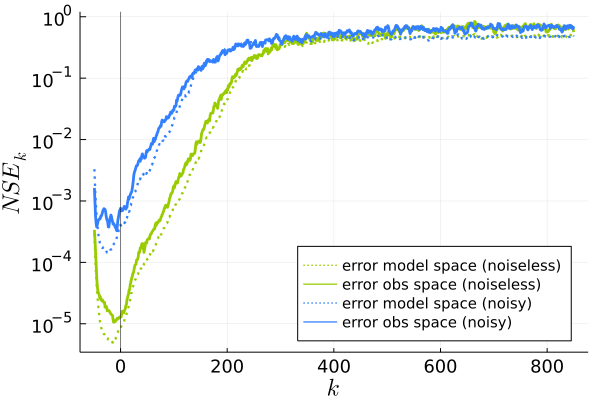}
	\caption{\textbf{Prediction error for the Lorenz system,} showing the median normalized squared error over $1000$ experiments for the observation space (solid lines) and the model space (dotted lines) for the case of noiseless (green) and noisy (blue) observations. The solid vertical line separates the assimilative regime ($k<0$) from the predictive regime ($k \geq 0$).}
	\label{fig:lorenz_error}
\end{figure}

We perform $1000$ independent experiments.  For each experiment we sample ground-truth microstates at random from the attractor of the system and generate time series of $T = 50$ observations with a sampling interval of $m = 2$ time steps per observation. We initialize the microstate for each time series following the steps presented 
in Section \ref{sec:summary}
using the parameters summarized in Table \ref{tab:params}. We present our results in Fig. \ref{fig:lorenz_error}, where we show the median assimilation ($k<0$) and prediction errors ($k \geq0$) of the noiseless and noisy time series for both the model and observation spaces (see Eqs. (\ref{eq:nse_obs})-(\ref{eq:nse_mod})).

From the assimilation side, our estimations get progressively better the closer they are to the present time at $k = 0$. This means that the longer the time series --and the more information we have from the \textit{past}--, the better the quality of the initialized microstate. Note that in the noisy case, the assimilative error plateaus near $10^{-3}$ in the observation space, marking the noise level of the observations. In contrast, the estimations keep getting better in the model space, indicating that even in the presence of noise, having more observations mitigates the probability of having dominating trajectories.

From the prediction side, we observe that the error diverges at essentially the same rate in both the noiseless and noisy cases. The main difference is that the error intercept at $k=0$ is higher in the noisy case, thus making the predictions saturate earlier than when the observations are noiseless. Specifically, we find prediction horizons of $k_{max} = 171$ and $k_{max}^{noisy} = 113$ steps, which correspond to $1.35 t_\lambda$ vs. $0.89 t_\lambda$ for the noiseless and noisy time series (see Fig \ref{fig:horizon_distribution} in Appendix \ref{app:interesting_figs} for an alternative approach on the prediction horizons). Moreover, we find that the prediction errors on the model and observation spaces are almost identical throughout the whole prediction window. Thus, measuring how the errors diverge in the observation space gives us a good proxy of the out-of-sample behavior of the latent microstates of the system. 

\begin{figure}[h!]
\centering
	\includegraphics[width=0.5\textwidth]{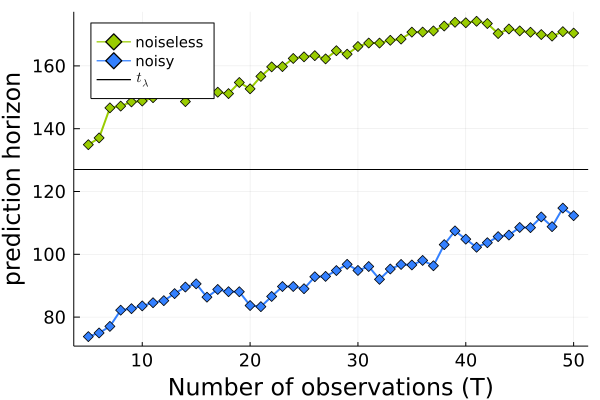}
	\caption{\textbf{Predictability vs. number of observations.} We show how the predictability horizon $k_{max}$ for the Lorenz system changes with the number of observations $T$ for noiseless (green) and noisy (blue) ensembles of time series. The horizontal black solid line indicates the Lyapunov 10-fold time $t_\lambda$.}
	\label{fig:lorenz_prediction_horizons}
\end{figure}

\begin{figure}[h!]
	\centering
	\includegraphics[width=0.49\textwidth]{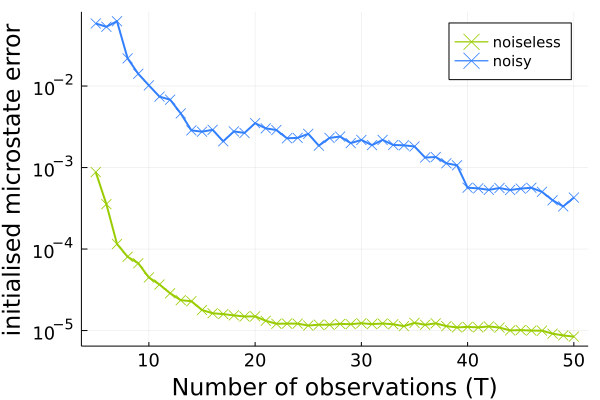}
	\caption{\textbf{Initialized microstate error for the Lorenz system.}  We show how the average discrepancy $\nse_0^{model}$ between the true present-time microstate $\bm{x}_0$ and the initialized microstate $\hat{\bm{x}}_0$ changes with the number of observations $T$ for noiseless (green) and noisy (blue) ensembles of time series.}
	\label{fig:lorenz_error_vs_measurements}
\end{figure}

In Figs. \ref{fig:lorenz_prediction_horizons} and \ref{fig:lorenz_error_vs_measurements} we analyze how the performance of our method depends on the length of the observation window, showing how the prediction horizons and the discrepancies between $\bm{x}_0$ and $\hat{\bm{x}}_0$ change with the number of observations. 

In Fig. \ref{fig:lorenz_prediction_horizons} we find that the prediction horizon increases linearly with the number of observations available, with a similar slope for both the noiseless and noisy cases. Not surprisingly, the noise affects the time horizon over which one can make an effective prediction. 

Additionally, we observe in Fig. \ref{fig:lorenz_error_vs_measurements} that the discrepancy 
decreases monotonically in both the noiseless and noisy cases. For the noiseless case, we observe a higher than exponential decrease in the discrepancy that ranges from $\nse_0^{model} \sim 10^{-3}$ for $T = 5$ to $\nse_0^{model} \sim 10^{-5}$ for $T=50$ observations. While the change in discrepancy is less pronounced for the noisy time series, it decreases 2 orders of magnitude with $\nse_0^{model} \sim 10^{-1.5}$ for $T = 5$ to $\nse_0^{model} \sim 10^{-3.5}$ for $T=50$ observations.

The Lorenz equations are a low-dimensional system ($N_x = 3$) with a low dimensional attractor of dimension $N_{\mathcal{M}} = 2.06$. The number of observations in our experiments can be much larger than the dimension of the system.  When combined with the fact that this system does not have any severe degeneracies,  (see Appendix \ref{app:nonautonomous_linear}), we recover the ground-truth microstate precisely with only a handful of noiseless observations. Every additional observation is, in theory, redundant for finding $\bm{x}_{-T}$ but, in practice, measurements can have several sources of error such as observational noise or the finite precision of numerical integration methods. Each additional observation thus further averages out these errors, which probably why $k_{max}$ gets better proportionally to $T$. 

\subsection{High-dimensional System: Mackey-Glass}
\label{sec:mackeyglass}

The Mackey-Glass system \cite{mackey1977oscillation} describes the dynamics of the density of cells in the blood with the following delayed differential equation
\begin{equation}
	\dot{x} = \mathscr{F}\left( x, x_{t_d} \right) = \frac{a x_{t_d}}{1 + x_{t_d}^c} - b x.
	\label{eq:mackey_glass_dde} 
\end{equation}
The state $x_{t_d} = x(t-t_d)$ is the density of cells delayed by $t_d$ time units and $a$, $b$, and $c$ are parameters. It exhibits chaotic dynamics for $t_d > 16.8$ with $a = 0.2, \ b = 0.1$ and $c = 10$ \cite{farmer1982chaotic}. In terms of blood cell density, the chaotic regime represents a pathological behavior.

The evolution of Eq. (\ref{eq:mackey_glass_dde}) relies on knowing the state of $x$ in the continuous interval $[t-t_d, t]$, making its state space infinite-dimensional. However, we can approximate such state by taking $N_x$ samples at intervals of length $\Delta t = t_d/N_x$ and constructing the $N_x$-dimensional microstate vector
\begin{align}
	\xk{k} &= \left( (\bm{x}_{k})_1, \cdots, (\bm{x}_{k})_{N_x} \right) \nonumber \\
	       &= \big( x( t_k - \underbrace{N_x \Delta t}_{t_d}), \dots, x(t_k - \Delta t), x(t_k) \big) ,
	\label{eq:discrete_x}
\end{align}
where $(\bm{x}_{i})_k = x(t - (N_x-i)\Delta t)$. 

Using this vector, we can obtain trajectories of the Mackey-Glass system with any numerical integrator, for which we use the Euler method with a fixed step size of $\Delta t$ for simplicity. We can thus recast this approximate system with the following $N_x$-dimensional deterministic mapping
\begin{align}
 \xk{k+1} =  \model{\xk{k}} = \left\lbrace
 \begin{array}{llll}
  (\bm{x}_t)_{N_x} + \Delta t \mathscr{F}\big( (\bm{x}_t)_{N_x}, (\bm{x}_t)_{1} \big) \\
  (\bm{x}_{t+1})_{1} + \Delta t \mathscr{F}\big( (\bm{x}_{t+1})_{1}, (\bm{x}_t)_{2}) \big) \\ 
  \vdots \\
  (\bm{x}_{t+1})_{N_x-1} + \Delta t \mathscr{F} \big( (\bm{x}_{t+1})_{N_x - 1}, (\bm{x}_{t})_{N_x} \big),
  \end{array} \right.
  \label{eq:mackey_glass_map}
\end{align} 
using $\mathscr{F}$ as defined in Eq. (\ref{eq:mackey_glass_dde}). The microstate-space dimension $N_x$ of this mapping is determined by $t_d/\Delta t$, for which we take $N_x = 50$ and $t_d = 25$ so that the system exhibits chaotic dynamics. Under these settings, the system's Lyapunovs 10-fold time is $t_\lambda = 230$ samples while its attractor has a Kaplan-Yorke dimension of $N_\mathcal{M} = 2.34$.

\begin{figure}[h!]
\centering
	\includegraphics[width=0.49\textwidth]{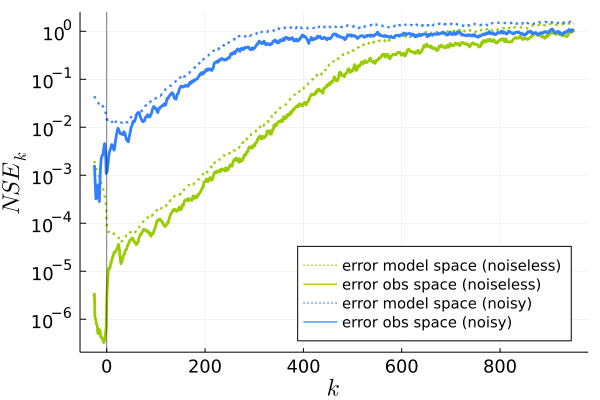}
    \caption{\textbf{Prediction error for the Mackey-Glass system.} We show the median normalized squared error over $1000$ experiments for the observation space (solid lines) and the model space (dotted lines) for the case of noiseless (green) and noisy (blue) observations. The solid vertical line separates the assimilative regime ($k<0$) from the predictive regime ($k \geq 0$).}
	\label{fig:mackey_error}
\end{figure}

As before, we perform $1000$ independent experiments in which, for each experiment, we sample ground-truth microstates at random from the attractor of the system and generate time series of $T=25$ observations with a sampling interval of $m = 2$ time steps per observation (see Table \ref{tab:params} for details). In contrast to the Lorenz system, we consider time series containing fewer data points than the dimension of the microstate space (in this case $T=25$ and $N_x=50$, respectively), making the problem under-determined. We present our results in Fig. \ref{fig:mackey_error}, where we show the median assimilation ($k<0$) and prediction errors ($k \geq0$) for the noiseless and noisy time series for both the model and observation spaces.

Unexpectedly, our method yields more accurate initializations for the Mackey-Glass system than the Lorenz system, even though both the attractor and the microstate space of the former have a higher dimension than the latter. However, if we compare the power spectra of the two systems (see Fig. \ref{fig:powerspectra} in Appendix \ref{app:interesting_figs}), we find that the Mackey-Glass system has a faster frequency decay and more frequency peaks than the Lorenz system, suggesting that the former is easier to initialize than the latter. For instance, the power amplitude for frequency $1/6$ is more than $1000$ higher in the Lorenz system than in the Mackey-Glass system. This further suggests that looking at the power spectra of the system is a better indicator of the initializability of a system than the dimension of its chaotic attractor.

From the prediction side, our results are qualitatively similar to what we saw for the Lorenz system. The error diverges at roughly the same rate in both the noiseless and noisy time series, with an error intercept at $k=0$ that is higher for the noisy experiments. Additionally, we find a very close correspondence between the errors in the model and observation spaces, which supports our claim that measuring the error in the observation space gives us a good proxy of the out-of-sample behavior of the latent microstate dynamics. 

In terms of their predictability horizons, we find that $k_{max} = 556$ and $k_{max}^{noisy} = 285$, corresponding to $2.4t_\lambda$ and $1.2t_\lambda$ for the noiseless and noisy ensembles respectively, (see Fig \ref{fig:horizon_distribution} in Appendix \ref{app:interesting_figs} for an alternative approach on the prediction horizons). We find it remarkable that with only $25$ observations of the system, we obtain predictions that stay accurate for significantly longer than the Lyapunov 10-fold time of the system. 

From the assimilation side, the microstate estimations get progressively better the closer they are to the present time, similar to what we observed in the Lorenz system. However, the assimilation error is significantly lower in the observation space than in the model space, suggesting, misleadingly, that the initialized microstate is much more accurate than the error we observe in the model space. Nonetheless, the errors in model and observation spaces converge to each other as soon as the prediction window starts, meaning that the error in the observation space is still an accurate proxy of the out-of-sample behavior in the model space.

Interestingly, the first few out-of-sample predictions in the model space have a lower discrepancy than in the observation space in both the noiseless and noisy cases. This happens, we believe, because the time span of the observations $\bm{y}$ is not long enough for the initialized microstate to converge onto the attractor. With only $T = 25$ data points of a $50$-dimensional chaotic system, we do not possess enough information to recover the present-time microstate precisely. 

\begin{figure}[h!]
\centering
	\includegraphics[width=0.5\textwidth]{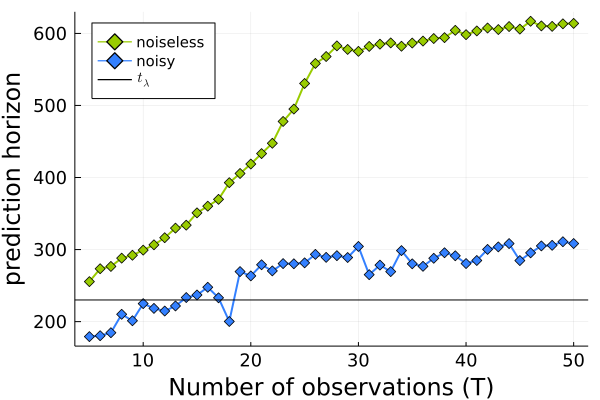}
	\caption{\textbf{Predictability horizon of the Mackey-Glass system.} We show how the predictability horizon $k_{max}$ changes with the number of observations $T$ for noiseless (green) and noisy (blue) ensembles of time series. The horizontal black solid line indicates the Lyapunov 10-fold time $t_\lambda$. For the noiseless case, we observe a critical transition on the behavior of $k_{max}$ for $T_c = 25$.}
	\label{fig:mackey_kmax}
\end{figure}

\begin{figure}[h!]
	\centering
	\includegraphics[width=0.49\textwidth]{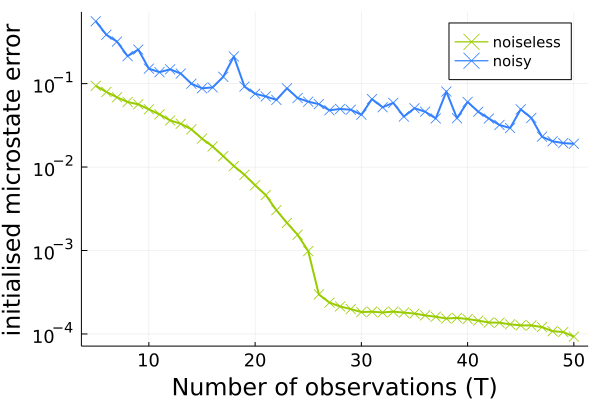}
	\caption{\textbf{Initialized model error for the Mackey-Glass system.} We show how $\nse_0^{model}$, the average discrepancy between the true present-time microstate $\bm{x}_0$ and the initialized microstate $\hat{\bm{x}}_0$, changes with the number of observations $T$ for noiseless (green) and noisy (blue) ensembles of time series. For the noiseless case, we observe a critical transition in the behavior of $\nse_0^{model}$ for $T_c = 25$.}
	\label{fig:mackey_error_vs_measurements}
\end{figure}

The previous discussion suggests that we need more observations to better initialize the system. To investigate this we perform a series of experiments in which we vary the number of data points of the observed time series and assess the quality of the predictions. Similar to what we did for the Lorenz system, we focus on the prediction horizon (see Fig \ref{fig:mackey_kmax}) and the discrepancy between the present-time and the initialized microstate (see Fig \ref{fig:mackey_error_vs_measurements}).

We find a contrasting behavior regarding the experiments between noisy and noiseless observations. When the observations are noisy, the prediction horizon increases linearly with the number of observations (see Fig. \ref{fig:mackey_kmax} blue), ranging from $k_{max} = 0.78 t_\lambda$ when $T=5$ to $k_{max} = 1.34 t_\lambda$ when $T=50$. We also find that, in general, the discrepancy between the initialized and ground-truth present microstate decreases monotonically with the number of observations (see Fig. \ref{fig:mackey_error_vs_measurements} bottom). These results are qualitatively similar to what we found for the Lorenz system: the more observations the better. 

When the observations are noiseless, we find a \textit{critical change of behavior} at roughly $T = 25$ observations. In Fig. \ref{fig:mackey_kmax} (green), we find that the prediction horizon rises superlinearly for $T < 25$, with $k_{max} = 1.11 t_\lambda$ for $T = 5$ to $k_{max} = 2.31 t_\lambda$ for $T = 25$. Afterwards, the prediction horizon grows linearly and with a marginal increase, getting to $k_{max} = 2.67 t_\lambda$ for $T = 50$. We note, however, that the increase in this linear regime is almost double of what we find in the noisy counterpart, with a slope of $\Delta k_{max} = 15.2$ steps per observation against $\Delta k_{max}^{noisy} = 8.3$, respectively. In parallel, we find that the discrepancy of the initialized microstate decreases abruptly for the time series of $ T \geq 25$ observations, as we show in Fig. \ref{fig:mackey_error_vs_measurements} (green). This sharp transition reflects that time series with more than $25$ observations have enough information to pin down the present-time latent microstate precisely, meaning that we possess enough data points to uniquely separate the mixing of the microstate generated by the observation operator $\mathcal{H}$ into its individual components.

\begin{figure}[h!]
\centering
	\includegraphics[width=0.5\textwidth]{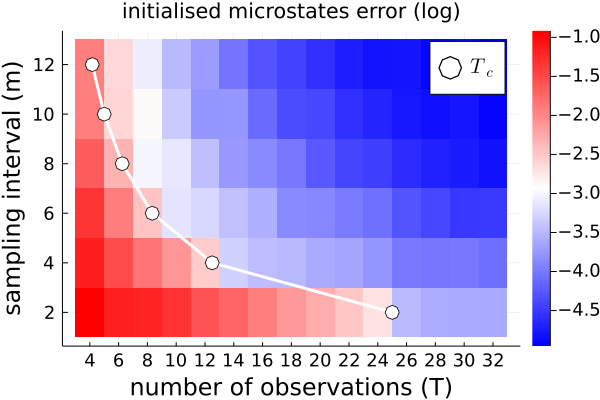}
	\caption{\textbf{Critical transition heatmap of the Mackey-Glass system:} In the $z$-axis, we show the (base-$10$ logarithm of the) initialized microstate discrepancy, $\nse_0^{model}$, as a function of the number of observations $T$ and the sampling interval $m$ for ensembles of noiseless time series. In white, we plot the $m = N_x/T$ curve for fixed $N_x = 50$. We find that the microstate discrepancies decrease abruptly before and after this curve.}
	\label{fig:mackey_error_heatmap}
\end{figure}

In short, having $25$ (or more) noiseless measurements of the system gives us enough information to precisely recover the present-time microstate, which is $50$-dimensional. Recall that we are considering the discrete map (\ref{eq:mackey_glass_map}) as the real system, so the only information we lose comes from either measuring the system with $\mathcal{H}$ or taking time series with a coarse-grained sampling frequency. Thus, for a fixed $\mathcal{H}$ and noiseless observations, recovering the initial microstates precisely should depend solely on how well the samples describe the latent trajectory of the system. Inspired by the Nyquist-Shannon sampling theorem \cite{shannon1949}, if we can establish a clear cutoff frequency on the power spectrum of the system, we could argue that if we observe the system with twice the frequency as the systems's cutoff frequency, then the observed signal would not lose any information with respect to the signal sampled for every update of map (\ref{eq:mackey_glass_map}). 

We thus claim that, if we can establish a clear cutoff frequency $f_c$ and the dynamical system does not suffer from severe degeneracies, we can precisely obtain the initial conditions of an $N_x$-dimensional (possibly) nonlinear system observed every $m$ updates with a scalar (possibly) nonlinear observation operator $\mathcal{H}$ if 1) $m^{-1} \geq 2f_c$ and 2) $T_c \geq N_x / m$. Having these conditions satisfied is equivalent to having an invertible observation-matrix $\bm{M}$ that determines the solution of the system $\bm{y} = \bm{M} \bm{x}_0$ exactly (see Appendix \ref{app:nonautonomous_linear} for a deeper development of this discussion).

We make further experiments to check the validity of our claim, where we measure the initialized microstate discrepancy when varying both the number of observations $T$ and the sampling interval $m$ while leaving the dimension of the system fixed to $N_x = 50$. If our claim about the Nyquist-Shannon theorem is a good approximation, we expect to find a $T_c \sim 1/m$ relation where, for $T \geq T_c$, the error in the initialized microstate becomes significantly lower than for $T < T_c$. In Fig \ref{fig:mackey_error_heatmap}, we show the results of these experiments, in which we observe such a change of regimes before and after the $T = N_x/m$ line, thus supporting the Nyquist-Shannon hypothesis.  

\section{Conclusions}\label{sec:conclusion}

Many natural and social processes can be accurately described by how their microstates interact and evolve into rich nontrivial dynamics with emerging properties. We often only possess aggregate noisy measurements of such processes, so it is of great interest to develop methods that let us extract information about the latent microstate dynamics from a given dataset. 

In this paper, we tackled the problem of initializing the latent microstate of a known (possibly) nonlinear dynamical system from aggregate (possibly) noisy and nonlinear observations of the system. We propose a three-step method to obtain such latent microstate that consists of 1) reducing the observational noise to mitigate possible dominating trajectories, 2) letting the system explore its attractor(s) and thus limiting the region in which we search for an optimal solution, and 3) minimizing the discrepancy between the simulated and real observations to obtain a refined estimation of the ground-truth microstate. We quantified the discrepancy between observations and simulations using a least-squares cost function in the observation space, similar to \cite{Pires1996, Carrassi2018}. We minimized the cost function using a plethora of gradient-based algorithms, for which we find that Adadelta \cite{zeiler2012adadelta} and Adam-oriented schemes  \cite{kingma2015adam, reddi2019convergence, yamada2018yamadam} perform the best.

We tested our method on two chaotic paradigmatic examples: the Lorenz and a high-dimensional approximation of the Mackey-Glass systems. We obtained initialized microstates that accurate fit the data, with out-of-sample predictions that outperformed the systems' Lyapunov 10-fold times, even when the observed time series were very short. We found that good predictions in the observations space always implied good predictions in the space of the microstates. We considered nonlinear observation operators that aggregate all the microstate component into a real number in all cases, with robust result with all the operators considered. Surprisingly, we obtained better results for the Mackey-Glass system, which has a higher-dimensional model space and higher-dimensional attractor but faster-decaying frequency spectrum than the Lorenz system, suggesting that the frequency spectrum gives us a better proxy of the initiability of the system than the observations to dimension of the model ratio. 

In most experiments, the quality of the initialized microstate was proportional to the number of data points of the observed time series. However, when the dimension of the system was higher than the number of observations and these observations were noiseless, the quality of the initialized microstate grew superlinearly. This superlinear regime transitions into the more common linear regime in a nontrivial manner, and we explored the conditions for such a transition. We claim that as long as we can establish a clear cutoff frequency of the observed data and this data meets the Nyquist-Shannon sampling theorem conditions with respect to the cutoff frequency, we can recover the ground-truth microstate precisely with fewer observations than the dimension of the system, thus marking the transition between regimes. This implies that if we possess a dataset where observations are sampled at an optimal rate such that we lose the least possible information of the underlying system, we can obtain high-quality initializations with just a handful of samples.

How well we can initialize a system depends on the amount of information the observed data contains and on the intrinsic features of the system. On the one hand, the amount of observational noise, the mixing that results from aggregating the system, the number of observations, and the data sampling rate contribute significantly in estimating microstates that may fit the data well but extrapolate poorly into the future. On the other hand, the dimension of the dynamical system, the frequency spectrum of its attractor(s), and how chaotic the system is, determines the window in which the predictions stay accurate.

This work gives us a conceptual framework to understand the interface between aggregate data and microscopic interactions. However, it is limited the case in which we know the model that perfectly specifies the system, as well as the observation operator and the characteristics of the noise in the observations. 

In future developments, we will explore how to deal with misspecified models and systems with stochastic behavior. We should include these new sources of uncertainty in the cost functions involved. With stochastic behavior alone, there are several new considerations for inferring individual agents' evolution in a system \cite{barde2012back}. Once we better understand what information low-dimensional observations gives us about a high-dimensional dynamical system, the natural next step is test our initialization method on real-world datasets.

\section{Statements and Declarations}

\textbf{Acknowledgments} The authors are grateful with Maria del Rio-Chanona, Marco Pangallo and Sylvain Barde for stimulating discussions, as well as the Oxford INET Complexity Economics group for valuable feedback. 
\\
\textbf{Code Availability} All the code described in the methodology is coded in Julia and is freely available at \url{https://github.com/blas-ko/LatentStateInitialization}.
\\
\textbf{Conflicts of Interest} The authors declare that they have no conflict of interest.
\\
\textbf{Funding} BK has received research support from the Conacyt-SENER: Sustentabilidad Energética scholarship.


\bibliographystyle{unsrt}
\bibliography{references}

\begin{thebibliography}{10}

\bibitem{Climate1999}
D~Rind.
\newblock Complexity and climate.
\newblock {\em Science}, 284:105--107, 1999.

\bibitem{ward2016dynamic}
Jonathan~A Ward, Andrew~J Evans, and Nicolas~S Malleson.
\newblock Dynamic calibration of agent-based models using data assimilation.
\newblock {\em Royal Society open science}, 3(4):150703, 2016.

\bibitem{Brain2009}
E~Bullmore and O~Sporns.
\newblock Complex brain networks: graph theoretical analysis of structural and
  functional system.
\newblock {\em Nature}, 10:186--198, 2009.

\bibitem{Markets2012}
JD~Farmer, M~Gallegati, C~Hommes, A~Kirman, P~Ormerod, S~Cincotti, A~Sanchez,
  and D~Helbing.
\newblock A complex systems approach to constructing better models for managing
  financial markets and the economy.
\newblock {\em European Physical Journal}, 214:295–324, 2012.

\bibitem{Networks2018}
Sabato Manfredi, Edmondo Di~Tucci, and Vito Latora.
\newblock Mobility and congestion in dynamical multilayer networks with finite
  storage capacity.
\newblock {\em Physical review letters}, 120(6):068301, 2018.

\bibitem{glaser2014complex}
Sarah~M Glaser, Michael~J Fogarty, Hui Liu, Irit Altman, Chih-Hao Hsieh, Les
  Kaufman, Alec~D MacCall, Andrew~A Rosenberg, Hao Ye, and George Sugihara.
\newblock Complex dynamics may limit prediction in marine fisheries.
\newblock {\em Fish and Fisheries}, 15(4):616--633, 2014.

\bibitem{grazzini2017bayesian}
Jakob Grazzini, Matteo~G Richiardi, and Mike Tsionas.
\newblock Bayesian estimation of agent-based models.
\newblock {\em Journal of Economic Dynamics and Control}, 77:26--47, 2017.

\bibitem{Farmer1980}
NH~Packard, JP~Crutchfield, JD~Farmer, and RS~Shaw.
\newblock Geometry from a time series.
\newblock {\em Physical Review Letters}, 45(9):712--716, 1980.

\bibitem{takens1981detecting}
Floris Takens.
\newblock Detecting strange attractors in turbulence.
\newblock In {\em Dynamical systems and turbulence, Warwick 1980}, pages
  366--381. Springer, 1981.

\bibitem{farmer1987predicting}
J~Doyne Farmer and John~J Sidorowich.
\newblock Predicting chaotic time series.
\newblock {\em Physical review letters}, 59(8):845, 1987.

\bibitem{sugihara2015free}
Hao Ye, R.~J. Beamish, S.~M. Glaser, S.~C.~H. Grant, Chih-hao Hsieh, L.~J.
  Richards, J.~T. Schnute, and G.~Sugihara.
\newblock Equation-free mechanistic ecosystem forecasting using empirical
  dynamic modeling.
\newblock {\em Proceedings of the National Academy of Sciences},
  116(41):E1569--E1576, 2015.

\bibitem{farmer1991noisered}
J~Doyne Farmer and John~J Sidorowich.
\newblock Optimal shadowing and noise reduction.
\newblock {\em Physica D: Nonlinear Phenomena}, 47(3):373--392, 1991.

\bibitem{Pires1996}
Carlos Pires, Robert Vautard, and Olivier Talagrand.
\newblock {On extending the limits of variational assimilation in nonlinear
  chaotic systems}, 1996.

\bibitem{Carrassi2018}
Alberto Carrassi, Marc Bocquet, Laurent Bertino, and Geir Evensen.
\newblock {Data assimilation in the geosciences: An overview of methods,
  issues, and perspectives}.
\newblock {\em Wiley Interdisciplinary Reviews: Climate Change}, 9(5):1--50,
  2018.

\bibitem{navon2009data}
Ionel~M Navon.
\newblock Data assimilation for numerical weather prediction: a review.
\newblock In {\em Data assimilation for atmospheric, oceanic and hydrologic
  applications}, pages 21--65. Springer, 2009.

\bibitem{platt2020comparison}
Donovan Platt.
\newblock A comparison of economic agent-based model calibration methods.
\newblock {\em Journal of Economic Dynamics and Control}, 113:103859, 2020.

\bibitem{evensen1997advanced}
Geir Evensen.
\newblock Advanced data assimilation for strongly nonlinear dynamics.
\newblock {\em Monthly weather review}, 125(6):1342--1354, 1997.

\bibitem{Judd2008}
Kevin Judd.
\newblock {Forecasting with imperfect models, dynamically constrained inverse
  problems, and gradient descent algorithms}.
\newblock {\em Physica D: Nonlinear Phenomena}, 237(2):216--232, 2008.

\bibitem{yang2010metaheuristics}
Xin-She Yang.
\newblock {\em Nature-inspired Metaheuristic Algorithms}.
\newblock Luniver Press, 2010.

\bibitem{ruder2016overview}
Sebastian Ruder.
\newblock An overview of gradient descent optimization algorithms.
\newblock {\em arXiv preprint arXiv:1609.04747}, 2016.

\bibitem{lorenz1963deterministic}
Edward~N Lorenz.
\newblock Deterministic nonperiodic flow.
\newblock {\em Journal of the atmospheric sciences}, 20(2):130--141, 1963.

\bibitem{Judd2007}
Kevin Judd.
\newblock {Failure of maximum likelihood methods for chaotic dynamical
  systems}.
\newblock {\em Physical Review E - Statistical, Nonlinear, and Soft Matter
  Physics}, 75(3):1--7, 2007.

\bibitem{Casdagli1991}
Martin Casdagli, Stephen Eubank, J.~Doyne Farmer, and John Gibson.
\newblock {State space reconstruction in the presence of noise}.
\newblock {\em Physica D: Nonlinear Phenomena}, 51(1-3):52--98, 1991.

\bibitem{guinon2007moving}
Jos{\'e}~Luis Gui{\~n}{\'o}n, Emma Ortega, Jos{\'e} Garc{\'\i}a-Ant{\'o}n, and
  Valent{\'\i}n P{\'e}rez-Herranz.
\newblock Moving average and savitzki-golay smoothing filters using mathcad.
\newblock {\em Papers ICEE}, 2007, 2007.

\bibitem{grassberger1993noise}
Peter Grassberger, Rainer Hegger, Holger Kantz, Carsten Schaffrath, and Thomas
  Schreiber.
\newblock On noise reduction methods for chaotic data.
\newblock {\em Chaos: An Interdisciplinary Journal of Nonlinear Science},
  3(2):127--141, 1993.

\bibitem{natrella2010nist}
Mary Natrella et~al.
\newblock e-handbook of statistical methods.
\newblock {\em NIST/SEMATECH}, 49, 2010.

\bibitem{anishchenko2002peculiarities}
Vadim~S Anishchenko, Tatjana~E Vadivasova, Andrey~S Kopeikin, J{\"u}rgen
  Kurths, and Galina~I Strelkova.
\newblock Peculiarities of the relaxation to an invariant probability measure
  of nonhyperbolic chaotic attractors in the presence of noise.
\newblock {\em Physical Review E}, 65(3):036206, 2002.

\bibitem{robbins1951stochastic}
Herbert Robbins and Sutton Monro.
\newblock A stochastic approximation method.
\newblock {\em The annals of mathematical statistics}, pages 400--407, 1951.

\bibitem{polyak1964some}
Boris~T Polyak.
\newblock Some methods of speeding up the convergence of iteration methods.
\newblock {\em Ussr computational mathematics and mathematical physics},
  4(5):1--17, 1964.

\bibitem{nesterov1983method}
Yurii~E Nesterov.
\newblock A method for solving the convex programming problem with convergence
  rate o (1/k\^{} 2).
\newblock In {\em Dokl. akad. nauk Sssr}, volume 269, pages 543--547, 1983.

\bibitem{duchi2011adaptive}
John Duchi, Elad Hazan, and Yoram Singer.
\newblock Adaptive subgradient methods for online learning and stochastic
  optimization.
\newblock {\em Journal of Machine Learning Research}, 12(Jul):2121--2159, 2011.

\bibitem{zeiler2012adadelta}
Matthew~D Zeiler.
\newblock Adadelta: an adaptive learning rate method.
\newblock {\em arXiv preprint arXiv:1212.5701}, 2012.

\bibitem{tieleman2012lecture}
Tijmen Tieleman and Geoffrey Hinton.
\newblock Lecture 6.5-rmsprop, coursera: Neural networks for machine learning.
\newblock {\em University of Toronto, Technical Report}, 2012.

\bibitem{kingma2015adam}
Diederik~P Kingma and Jimmy Ba.
\newblock Adam: A method for stochastic optimization.
\newblock {\em Proceedings of the 3rd International Conference on Learning
  Representations (ICLR)}, 2014.

\bibitem{reddi2019convergence}
Sashank~J Reddi, Satyen Kale, and Sanjiv Kumar.
\newblock On the convergence of adam and beyond.
\newblock {\em Proceedings of the 6th International Conference on Learning
  Representations (ICLR)}, 2019.

\bibitem{yamada2018yamadam}
Kazunori~D Yamada.
\newblock Yamadam: a hyperparameter-free gradient descent optimizer that
  incorporates unit correction and moment estimation.
\newblock {\em BioRxiv}, page 348557, 2018.

\bibitem{aurell1996growth}
Erik Aurell, Guido Boffetta, Andrea Crisanti, Giovanni Paladin, and Angelo
  Vulpiani.
\newblock Growth of noninfinitesimal perturbations in turbulence.
\newblock {\em Physical review letters}, 77(7):1262, 1996.

\bibitem{benettin1976kolmogorov}
Giancarlo Benettin, Luigi Galgani, and Jean-Marie Strelcyn.
\newblock Kolmogorov entropy and numerical experiments.
\newblock {\em Physical Review A}, 14(6):2338, 1976.

\bibitem{mackey1977oscillation}
Michael~C Mackey and Leon Glass.
\newblock Oscillation and chaos in physiological control systems.
\newblock {\em Science}, 197(4300):287--289, 1977.

\bibitem{farmer1982chaotic}
J~Doyne Farmer.
\newblock Chaotic attractors of an infinite-dimensional dynamical system.
\newblock {\em Physica D: Nonlinear Phenomena}, 4(3):366--393, 1982.

\bibitem{shannon1949}
C.~E. {Shannon}.
\newblock Communication in the presence of noise.
\newblock {\em Proceedings of the IRE}, 37(1):10--21, Jan 1949.

\bibitem{barde2012back}
Sylvain Barde.
\newblock Back to the future: economic rationality and maximum entropy
  prediction.
\newblock Technical report, School of Economics Discussion Papers, 2012.

\bibitem{tomas2010linear}
Mar{\'\i}a Tom{\'a}s-Rodr{\'\i}guez and Stephen~P Banks.
\newblock {\em Linear, time-varying approximations to nonlinear dynamical
  systems: with applications in control and optimization}, volume 400.
\newblock Springer Science \& Business Media, 2010.

\end{thebibliography}


\begin{appendices}

\section{Noiseless non-autonomous linear systems}\label{app:nonautonomous_linear}

In this section, we assume that the dynamical system (\ref{eq:ds}), the observations (\ref{eq:observations}), and the corresponding dynamical mapping (\ref{eq:dynamical_map}) are given by time-varying linear functions. We further assume that the dynamics are deterministic and the observations are noiseless so that $\bm{\xi}_k = 0$ and $\epsilon_k = 0$ for all $k$. Making these assumptions provides us with two advantages over arbitrary nonlinear dynamical systems: 1) linear systems are much easier to handle than nonlinear systems, and 2) nonlinear systems can be approximated by time-varying systems arbitrarily well if they are locally Lipchitz \cite{tomas2010linear}. 

Additionally, we find it convenient to slightly change the notation introduced in the main text. In the main text, we indexed the observations of the system so that the time stamps describe the observations $\bm{y}$ in the most natural way. Thus, we defined $t_k$ such that $y_k = y(t_k)$ is the $k$-th data point of the series. Consequently, we indexed the evolution of the microstates as $\bm{x}(t_{k+1}) = \bm{x}(t_k + m\Delta t) = \modelk{\bm{x}(t_k)}{}$, i.e., we needed to update the system $m$ times before sampling the next observation. Here, we index the passing of time in the time scale of the microstates, so that $\bm{x}(t_{k + 1}) = \bm{x}(t_k + \Delta t) = \model{\bm{x}(t_k)}$. Thus, we label the observed time series as $\bm{y} = (y_0, y_{m-1}, \dots, y_{mT -1}$) so that $y_{mk - 1}$ is the $k$-th data point of a time series of $T$ observations. This approach emphasizes that $\bm{y}$ is a coarse-grained sample of the underlying dynamics of the microstates.

Under the above considerations, we can recast Eqs. (\ref{eq:ds}-\ref{eq:observations}) respectively as follows
\begin{align}
    \bm{x}_{k + 1} &= \bm{F}_k \bm{x}_k = (\bm{F}_{k} \bm{F}_{k-1} \cdots \bm{F}_{0}) \xk{0}  \ , 
    \label{eq:linear_dynamics} \\
    y_k &= \bm{H}_k \xk{k} \ ,
    \label{eq:linear_observation}
\end{align}
where, at every time $t_k$, $\bm{F}_k$ is an ${N_x} \times {N_x}$ matrix representing a linear dynamical process and $\bm{H}_k$ is an $1 \times {N_x}$ matrix representing a linear observation operator. We assume that every element of $\bm{H}_k$ is non zero for every $k$. Note that under the current notation, the sequence $y_0, y_1, \dots, y_k, \dots$ represents the ground-truth dynamics under the (time-varying) observation operator $\bm{H}_k$, and only those indexes $k$ that are a multiple of $m$ are included in the observations $\bm{y}$.

From the RHS of Eq. (\ref{eq:linear_dynamics}), we see the explicit dependence of any observation $y_k$ from the initial conditions $\xk{0}$, so we can define an $1 \times N_x$ matrix\footnote{We define matrix $\bm{M}$ for time $t_0$ as $\bm{M}_0 := \bm{H}_0$.} 
\begin{equation*}
    \bm{M}_k := \bm{H}_k  \bm{F}_{k-1} \cdots \bm{F}_{0}    
\end{equation*}
that takes us from $\bm{x}_0$ to $y_k$ for any $k$. Thus, if we possess a time series of $T$ observations such that the last observations happens at time $t_{mT-1}$, then we can recast the microstate initialization problem as the following linear system of equations of $N_x$ variables and $mT$ equations
\begin{align*}
    \begin{aligned}
         y_0 &= \bm{M}_0 \xk{0} \\
         y_1 &= \bm{M}_1 \xk{0} \\
         &\vdots \\
         y_{mT-1} &= \bm{M}_{mT-1} \xk{0} \ , \\
    \end{aligned}
\end{align*}
or, more compactly,
\begin{equation}
    \bm{y}^* = \bm{M}^* \bm{x}_0 ,
    \label{eq:linear_extended_eq_system_matrix}
\end{equation}
where $\bm{y}^* = (y_0, y_1, \dots, y_{mT-1}) \in \mathbb{R}^{mT}$ is the \textit{extended} sequence of observations (it is extended in that it includes all the observations in $\bm{y} \in \mathbb{R}^T$ plus its intermediate, unobserved samples) and $\bm{M}^* = [\bm{M}_0 | \bm{M}_1 | \dots | \bm{M}_{mT-1}]$, the \textit{extended observation matrix} of size $mT \times N_x$.

We may solve system (\ref{eq:linear_extended_eq_system_matrix}) exactly whenever $\bm{M}$ is invertible. If the matrices $\bm{M}_k$ are non-singular, then $\bm{M}^*$ becomes invertible when $m T = N_x$.

We do not possess $\bm{y}^*$ but the coarse-grained time series $\bm{y}$ of $T$ observations, where two consecutive samples are $m$ time steps apart. Thus, we may only express a reduced form of system (\ref{eq:linear_extended_eq_system_matrix}) with $T$ equations and $N_x$ variables as
\begin{equation}
    \bm{y} = \bm{M} \bm{x}_0,
    \label{eq:linear_eq_system_matrix}
\end{equation}
with $\bm{M} = [\bm{M}_0 | \bm{M}_{m-1} | \dots | \bm{M}_{mT -1} ]$ of size $T \times N_x$. The system (\ref{eq:linear_eq_system_matrix}) spans the same time interval as system (\ref{eq:linear_extended_eq_system_matrix}), but it has $m$ less equations, so it is underdetermined and might have several solutions. 

Under the conditions we describe in what follows, we may obtain a solution $\bm{x}$ of the reduced system (\ref{eq:linear_eq_system_matrix}) that is also a solution of the extended system (\ref{eq:linear_extended_eq_system_matrix}). If $\bm{y}$ consists of $T = N_x /m$ (or more) data points, the solution $\bm{x}$ is unique and equal to the ground-truth microstate $\bm{x}_0$. More specifically, if the sampling frequency $(m \Delta t)^{-1}$ is higher than twice the cutoff frequency of the spectrum of the system and the matrices $\bm{M}_k$ are non singular for $k \in \{0, m-1, \dots, mT-1\}$, then the time series $\bm{y}$ determines the extended series $\bm{y}^*$ uniquely, and therefore any solution $\bm{x}$ that solves (\ref{eq:linear_eq_system_matrix}) also solves (\ref{eq:linear_extended_eq_system_matrix}). The above conditions establish the necessary and sufficient conditions for the Nyquist-Shannon sampling theorem \cite{shannon1949} to be true, so our result is a direct application of the theorem. 

Note that the power spectra of the systems considered in this work (and most chaotic systems) exhibit a power-law decay on their power spectrum, so there is not a well-defined cutoff frequency on neither the Mackey-Glass nor the Lorenz systems. Nevertheless, if their power spectrum decays fast enough, we should be able to define an effective cutoff frequency for the previous arguments to be a good approximation. We will investigate the validity of our arguments in future iterations of this work.

In this paper, we chose non-degenerate dynamical systems that are Lipschitz continuous and chose an observation operator (see Eq. (\ref{eq:observation_operator})) that is bijective to the linear operator $\bm{H} = [1, \dots, 1]^\dagger$. Therefore, the results of this Appendix apply in all our systems whenever the observations are noiseless. In particular, the results of this Appendix explain the critical transition we show in Figs. \ref{fig:mackey_kmax} and \ref{fig:mackey_error_vs_measurements} for $T = N_x / m = 25$.


\section{Comparison of several nonlinear observation operators}\label{app:observation_operators}

We assess the robustness of our initialization method using different (nonlinear) observation operators. The way we aggregate the microstates affects how much information we retain about the latent dynamics, so we consider the three following operators each with different levels of coupling between the microstate components
\begin{align}
    \mathcal{H}_1(\bm{x}) &= \sqrt[3]{ S_{\bm{x}} },
    \label{eq:observation_operator1} \\
    \mathcal{H}_2(\bm{x}) &= sign\left( P_{\bm{x}} \right) \sqrt[N_x]{ \lvert P_{\bm{x}} \rvert },
    \label{eq:observation_operator2} \\
    \mathcal{H}_3(\bm{x}) &= sign\left( C_{\bm{x}} \right) \sqrt{\lvert  C_{\bm{x}}  \rvert },
    \label{eq:observation_operator3}
\end{align}
where 
\begin{align}
    S_{\bm{x}} &= \sum_{i=1}^{N_x} \bm{x}_i^3, \\
    P_{\bm{x}} &= \prod_{i=1}^{N_x} \bm{x}_i, \\
    C_{\bm{x}} &= \sum_{i<j} \bm{x}_i \bm{x}_j.
\end{align}

Our rationale for choosing these operators goes as follows. We point out that Eq. (\ref{eq:observation_operator1}) is the same as Eq. (\ref{eq:observation_operator}): it's the cubic root of sum of cubes of the components of the microstates. This operator has no coupling between the microstate variables. As we mentioned before, Eq(\ref{eq:observation_operator1}) is bijective to any non-degenerate linear operator. For Eq. (\ref{eq:observation_operator2}), we take an operator that couples all the microstate components by multiplying them. If all the microstate components are positive --i.e. if $\bm{x}_i > 0$ for all $i$--, then we can take the logarithm of $\mathcal{H}_2$ and recover a non-degenerate linear operator. If any of the components is non-positive, then we cannot transform $\mathcal{H}_2$ to a linear operator, making the decoupling impossible. Finally, Eq. (\ref{eq:observation_operator3}) is the sum of pairwise couplings between the microstate components. In this case, there is no smooth transformation between $\mathcal{H}_3$ and a linear operator, so we can make no further simplification of $\mathcal{H}_3$. We consider the cubic-, $N_x$th-, and square-root of $S_{\bm{x}}$, $P_{\bm{x}}$, and $C_{\bm{x}}$ respectively so that the physical units of the observations are the same as the units of the microstates.

\begin{figure*}[h!]
\centering
	\begin{subfigure}[b]{\linewidth}
	\includegraphics[width=0.32\linewidth]{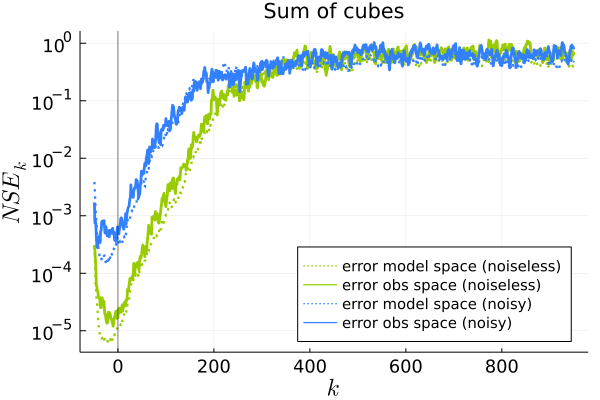}
	\includegraphics[width=0.33\linewidth]{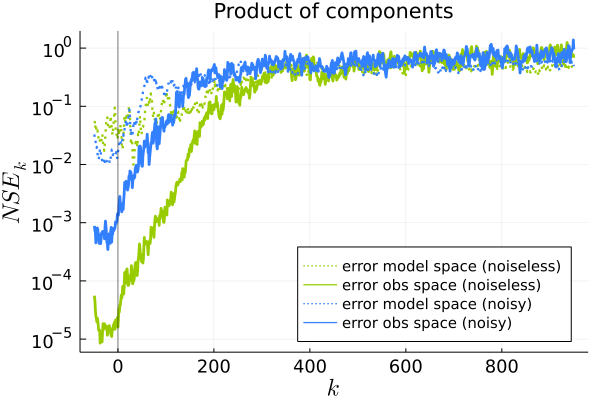}
	\includegraphics[width=0.33\linewidth]{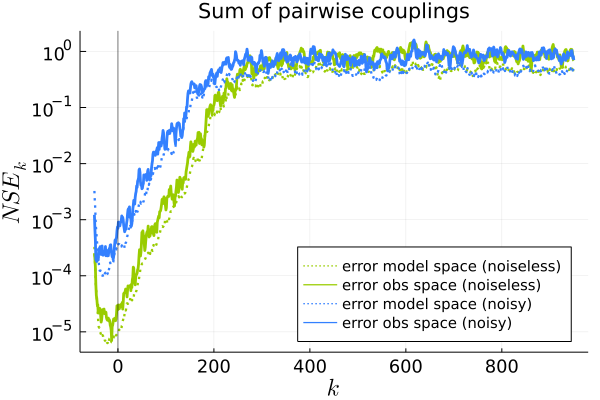}
	\caption{ Lorenz }	
	\end{subfigure}
	\begin{subfigure}[b]{\linewidth}
	\includegraphics[width=0.32\linewidth]{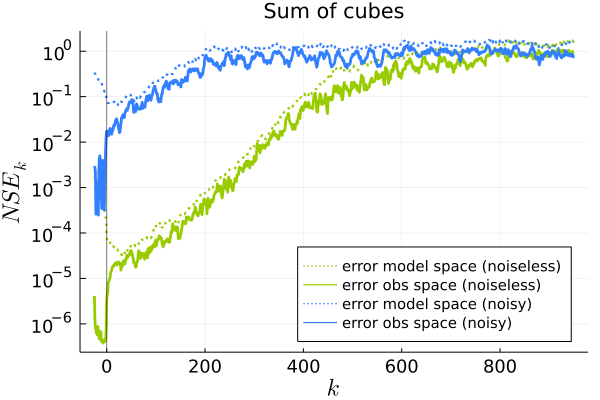}
	\includegraphics[width=0.33\linewidth]{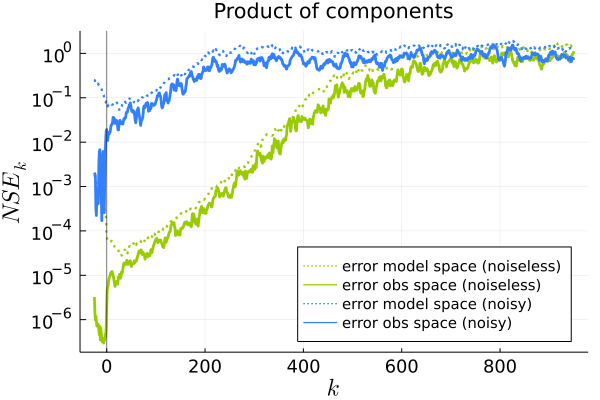}
	\includegraphics[width=0.33\linewidth]{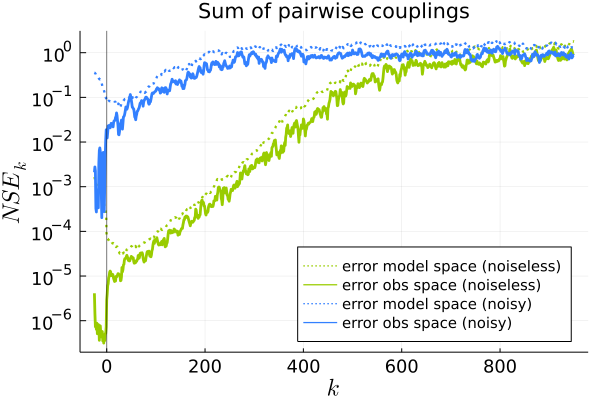}
	\caption{ Mackey-Glass }	
	\end{subfigure}
	\caption{\textbf{Prediction error for different observation operators.} We show the median normalized squared error over $200$ experiments for the observation space (solid lines) and the model space (dotted lines) for the case of noiseless (green) and noisy (blue) observations for \textit{a)} the Lorenz system and \textit{b)} the Mackey-Glass system. From left to right, we show the behaviour for operators representing the sum of cubes (see Eq. (\ref{eq:observation_operator1})), the product between microstate components (see Eq. (\ref{eq:observation_operator2})), and the sum of pairwise couplings (see Eq. (\ref{eq:observation_operator3})). We take all parameters as in Table \ref{tab:params}.}
	\label{fig:observation_operator_comparison}
\end{figure*}

In Fig. \ref{fig:observation_operator_comparison}, we show the median assimilation ($k<0$) and prediction ($k\geq0$) errors over $200$ runs for the Lorenz and the Mackey-Glass systems for each of the observation operators described above. We took the same set of $200$ initial conditions and noise seeds for each operator to make the results comparable.

In almost every case, we observe that the results for each system are almost identical regardless of the operator we use, both in the observation and the model spaces. This behavior suggests that our method is robust to the observation operator: if the observation convey enough information about the latent dynamics of the system, our method will initialize a microstate with good prediction power.

\begin{figure}[h!]
\centering
	\includegraphics[width=0.5\textwidth]{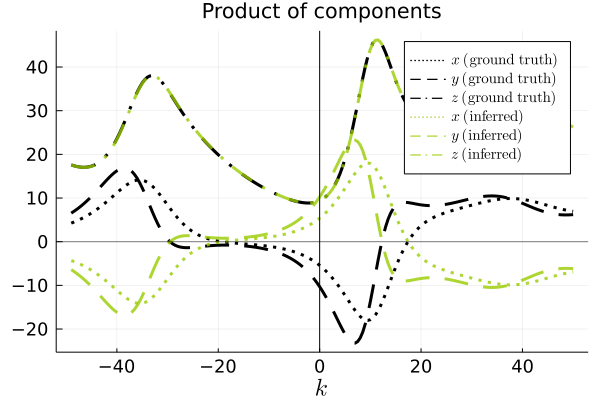}
	\caption{\textbf{Dynamics of the individual microstate components of the Lorenz system} for an example run of the ground truth (black) and the initialized (green) microstate when the observation are taken using the operator of Eq. (\ref{eq:observation_operator2}) --i.e., the operator that multiplies all the microstate components-- in a noiseless scenario. The $y$ and $z$ components are symmetric with respect to the $0$ line. Similar to Fig. \ref{fig:observation_operator_comparison}, $k <0$ denotes assimilation times and $k \geq 0$ prediction times.}
	\label{fig:component_symmetry_lorenz}
\end{figure}

However, we observe an anomaly on the model space error for the Lorenz system: the model space error corresponding to $\mathcal{H}_2$ --the product of components-- is significantly higher than i) its observation space error, and ii) the model error corresponding to $\mathcal{H}_1$ and $\mathcal{H}_3$. Recall that the Lorenz system is symmetric around its $x$-axis, so the operator $\mathcal{H}_2$ cannot resolve if the product $\prod_i \bm{x}_i$ corresponds to $(x,y,z)$ or to $(x,-y,-z)$, regardless of the number of measurements we have about the system. Thus, our method sometimes initializes the ground truth microstate and the other times it initializes its reflection about the $x$-axis, as we show in Fig. \ref{fig:component_symmetry_lorenz}.

Summarizing, the feasible set of solutions for the 
Lorenz system observed through $\mathcal{H}_2$ include both $(x,y,z)$ and $(x,-y,-z)$ which, incidentally, result in the same observation-space dynamics for time windows of arbitrary size. In this case, the initialization is more difficult because some microstates are indistinguishable. For all the other observation operators, the only feasible solution is the ground-truth microstate, hence their results are almost identical in \ref{fig:observation_operator_comparison}. These results suggest, that it is possible to improve the initialization procedure in systems with symmetries. It would only be necessary to add an extra step to the initialization procedure. This extra step would involve the rotation of the microstate found across all the symmetries of the system. We would then continue refining in parallel all those different symmetric microstates. The global solution would be then the one with the lowest cost function. We leave a deep analysis of this scenario for a future study.


\pagebreak
\section{Initialization method performance \& system features}\label{app:interesting_figs}

\begin{figure}[h!]
\centering
	\includegraphics[width=0.49\textwidth]{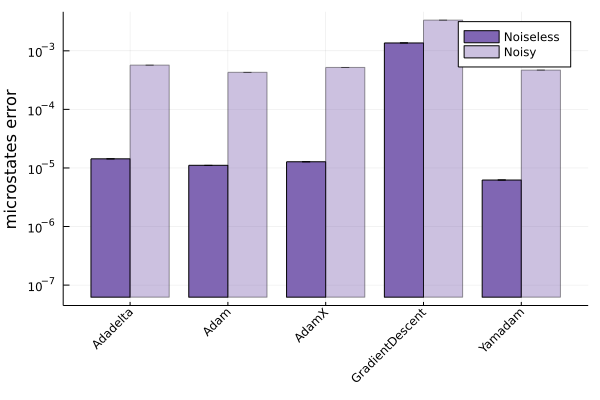}
	\caption{\textbf{Gradient-based optimizers performance on the Lorenz system.} We measure the performance of the gradient-based optimizers (described in section \ref{sec:refine}) with $\nse_0^{model}$, the average discrepancy between the present-time microstate $\bm{x}_0$ and the initialized microstate $\hat{\bm{x}}_0$, for noiseless (dark purple) and the noisy (light purple) time series. We show only the four best performing optimizers --namely Adadelta, Adam, AdamX, and YamAdam-- as well as Stochastic Gradient Descent, which serves as our benchmark.}
	\label{fig:lorenz_optimizers_comparison}
\end{figure}

\begin{figure}[h!]
\centering
	\includegraphics[width=0.49\textwidth]{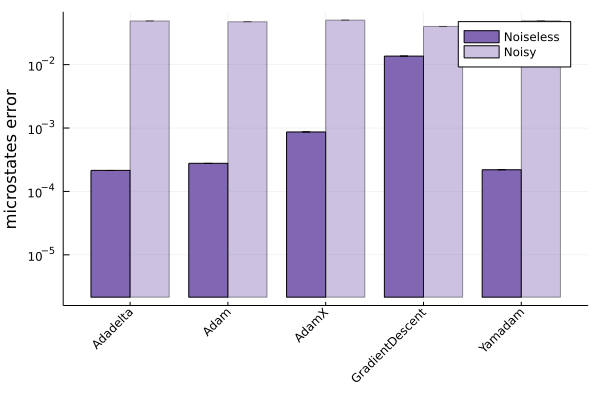}
	\caption{\textbf{Gradient-based optimizers performance on the Mackey-Glass system.} See description in Fig. \ref{fig:lorenz_optimizers_comparison} for details. The best performing optimizers were the same as with the Lorenz system.}
	\label{fig:mackey_optimizers_comparison}
\end{figure}

\begin{figure}[h!]
\centering
	\begin{subfigure}[b]{0.49\textwidth}
	\includegraphics[width=\textwidth]{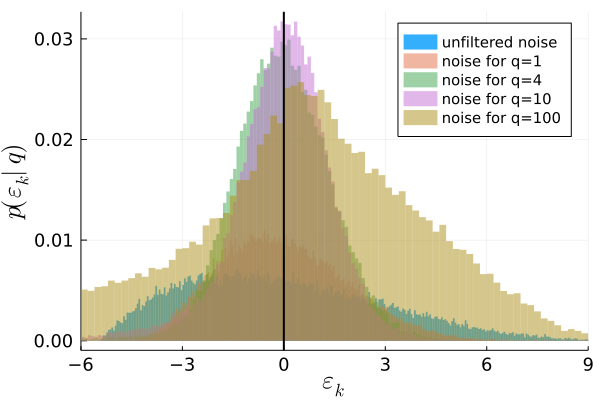}
	\caption{ Lorenz }	
	\end{subfigure}
	\begin{subfigure}[b]{0.49\textwidth}	
	\includegraphics[width=\textwidth]{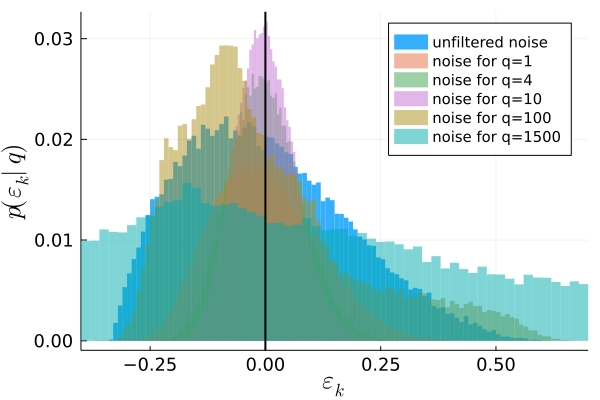}
	\caption{ Mackey-Glass }	
	\end{subfigure}
	\caption{\textbf{Filtered noise distributions} Starting from a simulated observational noise described by a left-skewed Beta distribution (see ``unfiltered noise'' in the legend) of very long time series ($T = 50000$ samples) of the \textit{a)} Lorenz and \textit{b)} Mackey-Glass systems, we plot the noise distribution of the unfiltered noise (solid blue) as well as the noise distribution we obtain after smoothing the corrupted signal with the LMPA filter for filters of increasing strength $q$. See section \ref{sec:preprocess} for details on the filter. For small $q$, the resulting distribution looks like a Gaussian distribution while for $q \gg 1$, the filter takes a significant of the signal, resulting in exhotic-shaped distributions.}
	\label{fig:error_distributions}
\end{figure}

\begin{figure}[h!]
\centering
	\begin{subfigure}[b]{0.49\textwidth}
	\includegraphics[width=\textwidth]{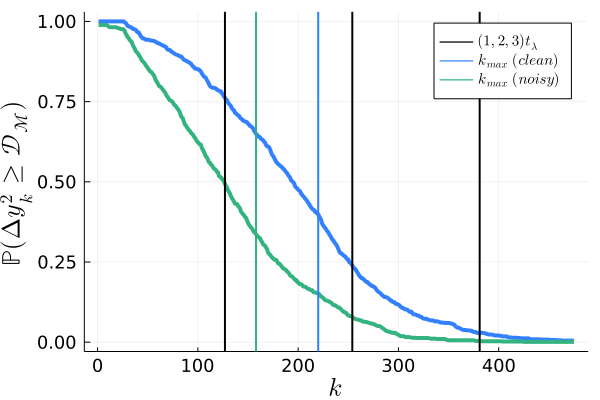}
	\caption{ Lorenz }	
	\end{subfigure}
	\begin{subfigure}[b]{0.49\textwidth}	
	\includegraphics[width=\textwidth]{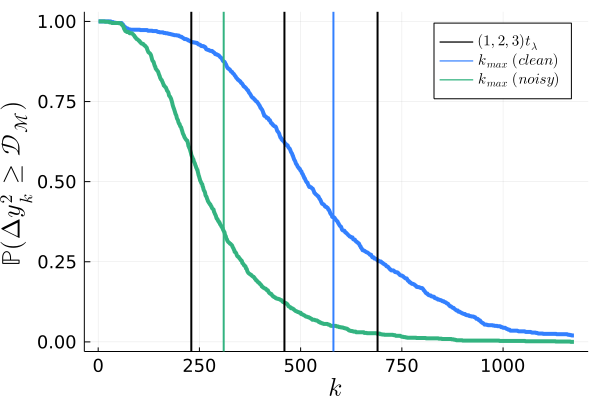}
	\caption{ Mackey-Glass }	
	\end{subfigure}
	\caption{\textbf{Cumulative probability of divergence} for \textit{a)} the Lorenz system and \textit{b)} the Mackey-Glass system. We show the fraction of trajectories --out of $1000$-- for which $NSE_k^{obs} \geq 2$ as a function of the prediction step $k$. This approach generalizes our definition of prediction horizon of Eq. (\ref{eq:k_max}) into a distribution-like quantity. In solid vertical lines, we show the Lyapunov 10-fold time of the system as well as the prediction horizons $k_{max}$ for the noiseless and noisy cases.}
	\label{fig:horizon_distribution}
\end{figure}

\begin{figure}[h!]
\centering
	\begin{subfigure}[b]{0.49\textwidth}
	\includegraphics[width=\textwidth]{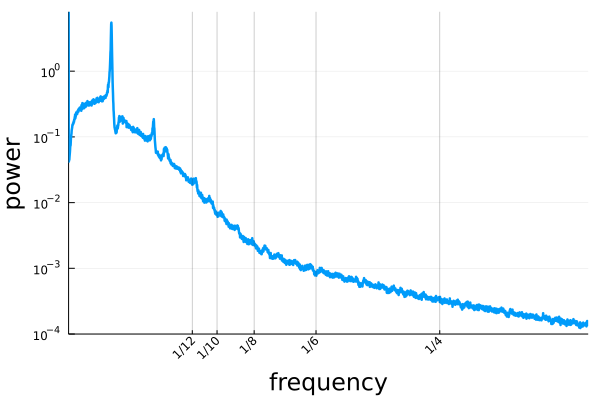}
	\caption{ Lorenz }	
	\end{subfigure}
	\begin{subfigure}[b]{0.49\textwidth}	
	\includegraphics[width=\textwidth]{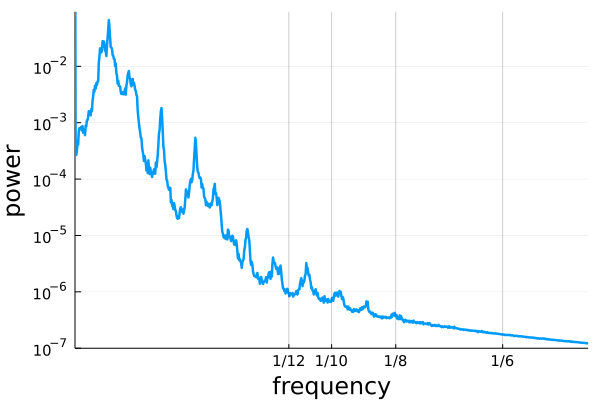}
	\caption{ Mackey-Glass }	
	\end{subfigure}
	\caption{\textbf{Normalized power spectra} for \textit{a)} the Lorenz system and \textit{b)} the Mackey-Glass system. The plot shows the average power spectra over $100$ random in-attractor initial conditions with trajectories of $2^{12}$ points. We note that the Lorenz system has a clear power-law frequency decay with only a few low frequency peaks. While the Mackey-Glass system exhibits a non-vanishing spectrum characteristic of chaotic systems, it has more defined frequency peaks and decays much faster than the Lorenz system. Thus, we expect for the Mackey-Glass system to be easier to initialize in general.}
	\label{fig:powerspectra}
\end{figure}

\begin{figure}[h!]
\centering
	\begin{subfigure}[b]{0.49\textwidth}
	\includegraphics[width=\textwidth]{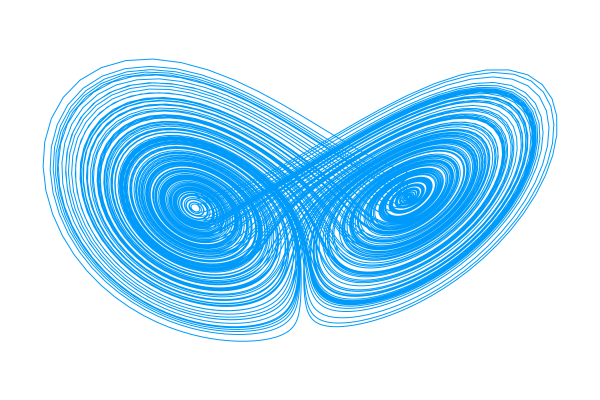}
	\caption{ Lorenz }	
	\end{subfigure}
	\begin{subfigure}[b]{0.49\textwidth}	
	\includegraphics[width=\textwidth]{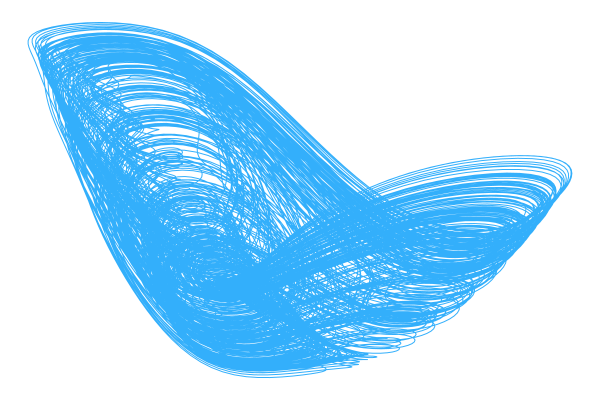}
	\caption{ Mackey-Glass }	
	\end{subfigure}
	\caption{\textbf{Chaotic attractors} for \textit{a)} the state space portrait of the Lorenz system, where the axis show $x, y$ and $z$, respectively, and \textit{b)} the state space portrait of the Mackey-Glass system, where we show $x(t)$ against $x(t - t_d)$. Each attractor consists of a \textit{long} trajectory of $15 000$ points coming from a random in-attractor initial condition.}
	\label{fig:attractors}
\end{figure}


\begin{figure*}[h!]
\centering
	\begin{subfigure}[b]{0.49\linewidth}
	\includegraphics[width=\linewidth]{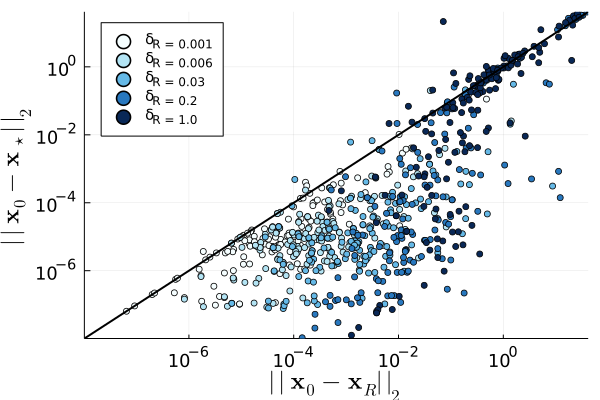}
	\includegraphics[width=\linewidth]{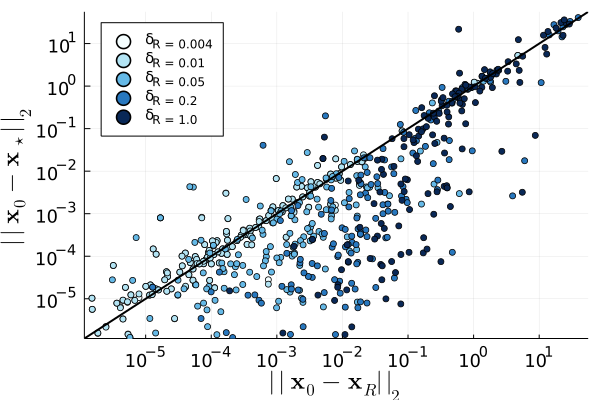}
	\caption{ Lorenz }	
	\end{subfigure}
	\begin{subfigure}[b]{0.5\linewidth}
	\includegraphics[width=\linewidth]{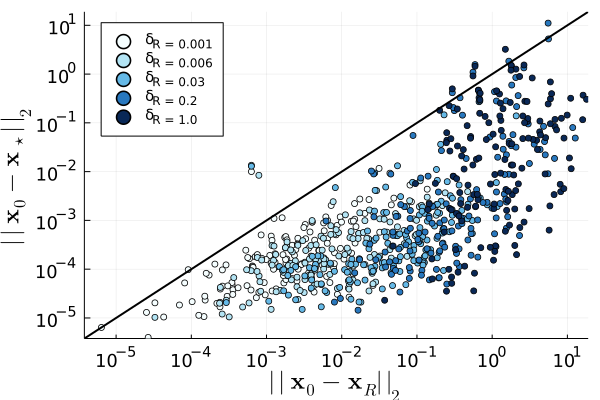}
	\includegraphics[width=\linewidth]{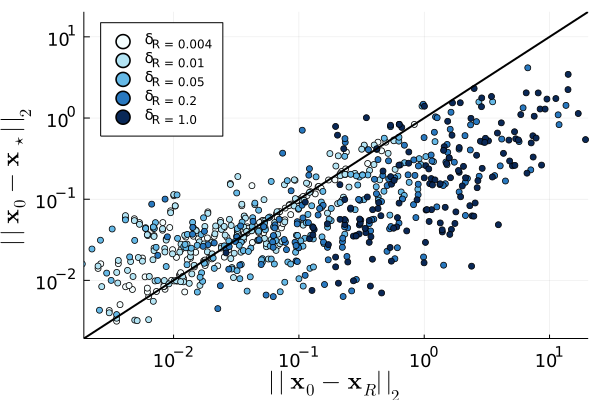}
	\caption{ Mackey-Glass }	
	\end{subfigure}
	\caption{\textbf{Initialization performance for different bounding-stage parameters} We present, on the $x$ axis, the discrepancy of $\hat{\bm{x}}_0^R$ against, on the $y$ axis, the discrepancy of the initialized microstates $\hat{\bm{x}}_0$ for increasing levels of $\delta_R$ (see legend) on \textit{a)} the Lorenz system and \textit{b)} the Mackey-Glass system. The microstate $\hat{\bm{x}}_0^R$ is the \textit{rough} estimation of $\bm{x}_0$ after Eq. (\ref{eq:deltaR}) is zatisfied. Values below the identity mean that $\hat{\bm{x}}_0^R$ improves refining it as described in section \ref{sec:refine}. Values on the diagonal mean that the initialized microstates do not get any better by applying refinement methods.}
	\label{fig:effects_of_bounding}
\end{figure*}

\end{appendices}

\end{document}